\documentclass[11pt,british,a4paper,reqno]{amsart}
    \usepackage[T1]{fontenc}
\usepackage[utf8]{inputenc}
	\usepackage[british]{babel}
    \DeclareMathSizes{12}{12}{7}{6}
	\usepackage{amsmath}
	\usepackage{amssymb}
	\usepackage{amsthm}
	\usepackage{bbm}
	\usepackage{braket}
	\usepackage{calc}
	\usepackage{cancel}
	\usepackage{xcolor}
	\usepackage[shortlabels]{enumitem}
	\usepackage{fancyhdr}
	\usepackage{float}
	\usepackage{graphicx}
	\usepackage{ifsym}
	\usepackage{indentfirst}
	\usepackage{mathrsfs}
	\usepackage{mathtools}
	\usepackage{parskip}
	\usepackage{proof}
	\usepackage{qtree}
	\usepackage{scalerel}
	\usepackage[%
	]{setspace}
	\usepackage{tikz}
	\usepackage{tikz-cd}
	\usepackage[normalem]{ulem}

	\usepackage{wasysym}
	\usepackage{oplotsymbl} %\pentago
	
\usepackage[hidelinks]{hyperref}
\usepackage[capitalise]{cleveref}
		\usepackage{bookmark}
			\bookmarksetup{numbered,open}

	\usepackage[%
	]{setspace}
\usepackage{geometry}
\bookmarksetup{numbered,open}

	\usetikzlibrary{matrix,arrows,decorations.pathmorphing,positioning,decorations.pathreplacing,shapes}
	\tikzset{commutative diagrams/.cd, 
		mysymbol/.style = {start anchor=center, end anchor = center, draw = none}}
	\tikzcdset{every label/.append style = {font = \footnotesize}}
	\tikzcdset{scale cd/.style={every label/.append style={scale=#1},
    cells={nodes={scale=#1}}}}

\newcommand{\BE}{\mathbb{E}}
\newcommand{\bN}{\mathbb{N}}

\newcommand{\A}{\mathcal{A}}
\newcommand{\B}{\mathcal{B}}
\newcommand{\CC}{\mathcal{C}}
\newcommand{\E}{\mathcal{E}}
\newcommand{\F}{\mathcal{F}}

\newcommand{\M}{\mathcal{M}}
\newcommand{\cP}{\mathcal{P}}

\newcommand{\fs}{\mathfrak{s}}

\newcommand{\X}{\mathcal{X}}
\newcommand{\Y}{\mathcal{Y}}
\newcommand{\Z}{\mathcal{Z}}

\newcommand{\Sim}{\operatorname{\mathsf{Sim}}\nolimits}
\newcommand{\Atom}{\operatorname{\mathsf{Atom}}\nolimits}
\newcommand{\Ab}{\operatorname{\mathsf{Ab}}\nolimits}
\newcommand{\add}{\operatorname{\mathsf{add}}\nolimits}
\newcommand{\proj}{\operatorname{\mathsf{proj}}\nolimits}
\newcommand{\Ext}{\operatorname{Ext}\nolimits}
\newcommand{\Hom}{\operatorname{Hom}\nolimits}

\newcommand{\Iso}{\operatorname{Iso}\nolimits}

\newcommand{\Sym}[1]{\operatorname{Sym}\nolimits({#1})}

\newcommand{\bsm}{\begin{psmallmatrix}}
\newcommand{\esm}{\end{psmallmatrix}}

\newcommand{\op}{\mathrm{op}}

\newcommand{\iso}{\cong}
\newcommand{\niso}{\ncong}
\newcommand{\into}{\hookrightarrow}

\newcommand{\sse}{\subseteq}

\newcommand{\deff}{\coloneqq}
\renewcommand{\geq}{\geqslant}
\renewcommand{\leq}{\leqslant}
\renewcommand{\phi}{\varphi}
\newcommand{\eps}{\varepsilon}
\newcommand{\ring}{R}

\newcommand\restr[2]{{\left.\kern-\nulldelimiterspace#1
						\right|_{#2}}}
						
\newcommand{\id}[1]{\mathrm{id}_{#1}}
\newcommand{\lmod}[1]{{#1}\operatorname{--\,\mathsf{mod}}\nolimits}

\newcommand{\dExt}[1]{\operatorname{-Ext}\nolimits(#1)}

	\newcommand{\ul}[1]{\underline{#1}}

\newcommand{\overlap}[2]{%
  \leavevmode\begingroup
  \vphantom{#1#2}%
  \ooalign{\hfil$#1$\hfil\cr\hfil$#2$\hfil\cr}%
  \endgroup
}

\DeclarePairedDelimiter\abs{\lvert}{\rvert}%

    \makeatletter
        \let\oldabs\abs
        \def\abs{\@ifstar{\oldabs}{\oldabs*}}
    \makeatother

	\theoremstyle{plain}
		\newtheorem{thm}{Theorem}[section] % reset thm numbering for each section
			\crefname{thm}{Theorem}{Theorems}
		\newtheorem{prop}[thm]{Proposition}
			\AddToHook{env/prop/begin}{\crefalias{thm}{prop}}
			\crefname{prop}{Proposition}{Propositions}
		\newtheorem{lem}[thm]{Lemma}
			\AddToHook{env/lem/begin}{\crefalias{thm}{lem}}
			\crefname{lem}{Lemma}{Lemmas}
		\newtheorem{cor}[thm]{Corollary}
			\AddToHook{env/cor/begin}{\crefalias{thm}{cor}}
			\crefname{cor}{Corollary}{Corollaries}
		
		\newtheorem{thmx}{Theorem}
	\theoremstyle{definition}
		\newtheorem{definition}[thm]{Definition}
			\AddToHook{env/definition/begin}{\crefalias{thm}{definition}}
			\crefname{definition}{Definition}{Definitions}
			
		\newtheorem{example}[thm]{Example}
			\AddToHook{env/example/begin}{\crefalias{thm}{example}}
			\crefname{example}{Example}{Examples}

		\newtheorem{notation}[thm]{Notation}
			\AddToHook{env/notation/begin}{\crefalias{thm}{notation}}
			\crefname{notation}{Notation}{Notations}

		\newtheorem{setup}[thm]{Setup}
			\AddToHook{env/setup/begin}{\crefalias{thm}{setup}}
			\crefname{setup}{Setup}{Setups}

	\theoremstyle{remark}
		\newtheorem{remark}[thm]{Remark}
			\AddToHook{env/remark/begin}{\crefalias{thm}{remark}}
			\crefname{remark}{Remark}{Remarks}
		
	\numberwithin{figure}{section}
	\numberwithin{equation}{section}

\makeatletter

	\renewcommand{\andify}{%
		\nxandlist{\unskip, }{\unskip{} \@@and~}{\unskip \penalty-2 \space \@@and~}}
    
	\renewcommand\author@andify{%
  		\nxandlist {\unskip ,\penalty-1 \space\ignorespaces}%
		{\unskip {} \@@and~}%
		{\unskip \penalty-2 \space \@@and~}
	}

	    \newenvironment{acknowledgements}{%
	    \renewcommand\abstractname{\textbf{Acknowledgements}}
\global\setbox\abstractbox=\vtop \bgroup
\normalfont\Small
\list{}{\labelwidth\z@
\leftmargin3pc \rightmargin\leftmargin
\listparindent\normalparindent \itemindent\z@
\parsep\z@ \@plus\p@

}%
\item[\hskip\labelsep\abstractname.]%
}{%
\endlist\egroup
\ifx\@setabstract\relax \@setabstracta \fi
}

\def\@setaddresses{\par
    \nobreak \begingroup
    \setstretch{0.5} %% Changes spacing
    \footnotesize
    \def\author##1{\nobreak\addvspace\bigskipamount}%
    \def\\{\unskip, \ignorespaces}%
    \interlinepenalty\@M
    \def\address##1##2{\begingroup
        \par\addvspace\bigskipamount\indent
    \@ifnotempty{##1}{(\ignorespaces##1\unskip) }%
    {\scshape\ignorespaces##2}\par\endgroup}%
    \def\curraddr##1##2{\begingroup
    \@ifnotempty{##2}{\nobreak\indent\curraddrname
    \@ifnotempty{##1}{, \ignorespaces##1\unskip}\/:\space
    ##2\par}\endgroup}%
    \def\email##1##2{\begingroup
    \@ifnotempty{##2}{\nobreak\indent\emailaddrname
    \@ifnotempty{##1}{, \ignorespaces##1\unskip}\/:\space
    \ttfamily##2\par}\endgroup}%
    \def\urladdr##1##2{\begingroup
    \def~{\char'\~}%
    \@ifnotempty{##2}{\nobreak\indent\urladdrname
    \@ifnotempty{##1}{, \ignorespaces##1\unskip}\/:\space
    \ttfamily##2\par}\endgroup}%
    \addresses
    \endgroup
}

\let\oldtocsection=\tocsection
\let\oldtocsubsection=\tocsubsection

\renewcommand{\tocsection}[2]{\vspace*{-5pt}\hspace{0em}\oldtocsection{#1}{#2}}

    \renewcommand{\tocsubsection}[2]{\vspace*{-5pt}\hspace{21pt}\oldtocsubsection{#1}{#2}}
\makeatother

\title{Stratifying systems and Jordan-H{{\"{o}}}lder extriangulated categories}

\author[Br{\"{u}}stle]{Thomas Br{\"{u}}stle}
    \address{
        D{\'{e}}partment de math{\'{e}}matiques\\
        Universit{\'{e}} de Sherbrooke\\ 
        Sherbrooke, Qu{\'{e}}bec\\ %J1K 2R1 \\ 
        Canada
    }
    \email{thomas.brustle@usherbrooke.ca}

\author[Hassoun]{Souheila Hassoun}
    \email{souheila.hassoun@usherbrooke.ca}

\author[Shah]{Amit Shah}
	\address{
		Department of Mathematics\\
		Aarhus University\\
        Aarhus\\
		Denmark
	}
    \email{a.shah1728@gmail.com}
    
\author[Tattar]{Aran Tattar}
	\address{
		Mathematical Institut\\
        University of Cologne\\
        K{\"{o}}ln\\
        Germany
    }
    \email{aran.tattar@gmail.com}

\date{\today}



\dedicatory{Dedicated to the late Brian Parshall in honour of his contributions to representation theory.}

\begin{document}

\begin{abstract} 
Stratifying systems, which have been defined for module, triangulated and exact categories previously, were developed to produce examples of standardly stratified algebras. A stratifying system $\Phi$ is a finite set of objects satisfying some orthogonality conditions. One very interesting property is that the subcategory $\mathcal{F}(\Phi)$ of objects admitting a composition series-like filtration with factors in $\Phi$ has the Jordan-H{\"{o}}lder property on these filtrations. 

This article has two main aims. First, we introduce notions of subobjects, simple objects and composition series for an extriangulated category, in order to define a Jordan-H{\"{o}}lder extriangulated category. Moreover, we characterise Jordan-H{\"{o}}lder, length, weakly idempotent complete extriangulated categories in terms of the associated Grothendieck monoid and Grothendieck group. Second, we develop a theory of stratifying systems in extriangulated categories. We define projective stratifying systems and show that every stratifying system $\Phi$ in an extriangulated category is part of a minimal projective one $(\Phi,Q)$. We prove that $\mathcal{F}(\Phi)$ is a length, Jordan-H{\"{o}}lder extriangulated category when $(\Phi,Q)$ satisfies a left exactness condition.

We give several examples and answer a recent question of Enomoto--Saito in the negative. 
\end{abstract}

\maketitle

\vspace{-0.5cm}
{\setstretch{0.5}\tableofcontents}
\vspace{-1cm}

%%%%%%%%%%%%%%%%%%%%%%%%%%%%%%%%%%%%%%
%%%%%%%%%%%%%%%%%%%%%%%%%%%%%%%%%%%%%%

\section{Introduction}
\label{sec:introduction}

In groundbreaking work, Cline--Parshall--Scott \cite{ClineParshallScott-finite-dimensional-algebras-and-highest-weight-categories} proved that there is a one-to-one correspondence between 
equivalence classes of highest weight categories with finitely many simple objects
and 
Morita equivalence classes of quasi-hereditary algebras (see \cite[Thm.\ 3.6]{ClineParshallScott-finite-dimensional-algebras-and-highest-weight-categories}). 
Consequently, a strong connection between representation theory of Lie algebras and that of finite-dimensional algebras was established. 
Indeed, since then 
highest weight categories have attracted large interest in representation theory 
(see e.g.\ 
\cite{%
Doty-representation-theory-of-reductive-normal-algebraic-monoids,
Guay-projective-modules-in-the-category-O-for-the-cherednik-algebra,
BezrukavnikovEtingof-parabolic-induction-and-restriction-functors-for-rational-cherednik-algebras,
BellamyThiel-highest-weight-theory-for-finite-dimensional-graded-algebras-with-triangular-decomposition,
Coulembier-some-homological-properties-of-ind-completions-and-highest-weight-categories}). 
Another very interesting result in \cite{ClineParshallScott-finite-dimensional-algebras-and-highest-weight-categories} was the Brauer-Humphreys reciprocity theorem \cite[Thm.\ 3.11]{ClineParshallScott-finite-dimensional-algebras-and-highest-weight-categories}, which is an instance of the Jordan-H{\"{o}}lder phenomenon. 
Using the dimensions of some $\Hom$-spaces, it was shown that 
the number of times an  object $A_{\lambda}$
appears as a composition factor in a good filtration of a filtered object $I_{\lambda}$ is independent of the chosen good filtration. 
A similar result is obtained in \cite[Thm.\ 3.4.7]{ClineParshallScott-stratifying-endomorphism-algebras} using ranks of bilinear forms.
Only later was a more general proof of why these multiplicities are independent of a chosen filtration in module categories given by Erdmann--S\'{a}enz \cite{ES}. 
One of the main aims of our article is to extend this approach
and, in particular, study length and the Jordan-H{\"{o}}lder property in extriangulated categories.

The notion of a stratifying system was introduced in \cite{ES} for module categories, 
providing a way to generate examples of standardly stratified algebras in the sense of Xi 
\cite{Xi-standardly-stratified-algebras-and-cellular-algebras} (see also \cite{AgostonHappelLukacsUnger-standardly-stratified-algebras-and-tilting, ClineParshallScott-stratifying-endomorphism-algebras,  DlabRingel-the-module-theoretical-approach-to-quasi-hereditary-algebras}). 
These algebras generalise quasi-hereditary algebras, which were introduced by Scott  \cite{Scott-simulating-algebraic-geometry-with-algebra-I}. 
Similar systems for certain exact categories had already been considered in \cite[Lem.~ 2.3]{BH},
and
the notion of a stratifying system has also been developed by Mendoza--Santiago \cite{MS} and Santiago \cite{Sa} for triangulated and exact categories, respectively.
Building on this work and generalising several results, we study stratifying systems in the realm of extriangulated categories. 
Independently, Adachi--Tsukamoto \cite{adachi2022mixed} have recently investigated similar concepts in the extriangulated setting, namely those of  \emph{mixed standardisable sets}; see \cref{rem-definitioncomparison}. 

Extriangulated categories were defined by Nakaoka--Palu \cite{NP19}, and unify exact and triangulated frameworks. 
Suppose 
$(\A,\BE,\fs)$ 
is an \emph{artin} extriangulated category (see \cref{def:artin-extriangulated-category}). 
Recall that $(\A,\BE,\fs)$ being extriangulated means that $\A$ is an additive category and that 
$\BE\colon\A^{\op} \times \A \to \Ab$ is a biadditive functor to the category of abelian groups. 
Given an element $\delta\in\BE(C,A)$, 
the realisation $\fs$ assigns to $\delta$ an equivalence class 
$
\fs(\delta)
    = [\begin{tikzcd}[column sep=0.5cm]
    A \arrow{r}{a}& B \arrow{r}{b}& C
    \end{tikzcd}]
$
of a pair of composable morphisms in $\A$. 
This is denoted by the \emph{extriangle}
$
\begin{tikzcd}[column sep=0.5cm,cramped]
A \arrow{r}{a}& B \arrow{r}{b}& C \arrow[dashed]{r}{\delta}&{}.
\end{tikzcd}
$ 
In \cref{sec:extriangulated-categories} 
we recall the basics and a recent characterisation 
of extriangulated categories 
with a weakly idempotent complete underlying additive category; 
see \cref{prop:weaklyidempotent}.

We call a set 
$\Phi = \{\Phi_{i}\}_{i=1}^{n}$ 
of indecomposable objects of $\A$  
an \emph{$\BE$-stratifying system} 
if it satisfies some $\Hom$- and $\BE$-orthogonality; see \cref{def:stratifying-system}. 
Given such a collection of objects, it is natural to consider the full subcategory $\F(\Phi)$ of $\A$ consisting of the \emph{$\Phi$-filtered objects}, i.e.\ 
the objects $M\in\A$ admitting a \emph{$\Phi$-filtration} 
\begin{equation}
\label{eqn:Theta-filtration}
\begin{tikzcd}[column sep=0.35cm, row sep=0.5cm]
0\arrow{rr}&&0 = M_{0} \arrow{rr}\arrow{dl}&& M_{1} \arrow{rr}\arrow{dl}&& M_{2} \arrow{r}\arrow{dl}& \cdots \arrow{r}& M_{t-1} \arrow{rr}&& M_{t} = M \arrow{dl}\\
&0\arrow[dashed]{ul}{\xi_{0}}&&\Phi_{j_{1}}\arrow[dashed]{ul}{\xi_{1}}&&\Phi_{j_{2}}\arrow[dashed]{ul}{\xi_{2}} &&&&\Phi_{j_{t}}. \arrow[dashed]{ul}{\xi_{t}}&
\end{tikzcd}
\end{equation}
in $(\A,\BE,\fs)$, 
with 
$\Phi_{j_{i}} \in \Phi$ for $1\leq i\leq t$. 
The objects $\Phi_{j_{i}}$ in the filtration can sometimes be thought of as ``composition factors''; 
see 
\cref{cor:simple-isomorphic-to-X}. 
Studying filtered objects has a long tradition in module categories (e.g.\ 
\cite{%
Ringel-the-category-of-modules-with-good-filtrations-over-a-quasi-hereditary-algebra-has-almost-split-sequences,
KoenigKulshammerOvsienko-quasi-hereditary-algebras-exact-borel-subalgebras-A-infty-categories-and-boxes,
Conde-delta-filtrations-and-projective-resolutions-for-the-auslander-dlab-ringel-algebra}), 
and in particular one notices that a filtration by simple modules is the same as a composition series. 
Filtrations of objects also appear in 
triangulated \cite{Bridgeland-stability-conditions-on-triangulated-categories,MS}, 
exact \cite{Sa}, 
and 
extriangulated \cite{Zh19} categories. 
A categorical perspective on the module case has been taken by Bodzenta--Bondal \cite{BB} recently.

It is shown by Zhou \cite{Zh19} that $\F(\Phi)$ is extension-closed in $(\A,\BE,\fs)$ (see \cref{rem:comments-on-FTheta}\ref{part:F(X)-is-smallest-ext-closed-subcat-containing-X}). 
Since extension-closed subcategories of exact categories are again exact, we see that $\F(\Phi)$ would inherit an exact structure if $(\A,\BE,\fs)$  were exact. 
Thus, the subcategories of $\Phi$-filtered objects studied in \cite{Sa} are again exact. 
Furthermore, the systems considered in \cite{MS}  for triangulated categories and in \cite{adachi2022mixed} for the extriangulated case are required to satisfy a condition that implies that $\F(\Phi)$ is an exact category; 
see the discussion in \cref{rem-definitioncomparison}.
The reason for this somewhat strong restriction might be because extension-closed subcategories of triangulated categories are not necessarily triangulated. 
But this is where the utility of extriangulated category theory becomes especially apparent: extension-closed subcategories of triangulated categories
are always extriangulated (see \cref{example:extension-closed-subcat-of-extriangulated-is-extriangulated}). 
Thus, for an $\BE$-stratifying system $\Phi$ in the extriangulated category $(\A,\BE,\fs)$, we have that $\F(\Phi)$ inherits an extriangulated structure, which we denote by 
$(\F(\Phi), \restr{\BE}{\F(\Phi)},\restr{\fs}{\F(\Phi)})$. 
Moreover, our axioms for $\Phi$ do not imply $\F(\Phi)$ is exact in general; 
see \cref{sec:examples} for examples.

When $\Phi$ is an $\BE$-stratifying system, the elements of $\Phi = \{\Phi_{i}\}_{i=1}^{n}$ behave like simple objects in the subcategory $\F(\Phi)$. 
Motivated by concepts used in \cite{ES,LMMS,MS,Sa}, we introduce \emph{projective $\BE$-stratifying systems} which adds the data of projective-like objects in $\F(\Phi)$.
A projective $\BE$-stratifying system is a pair $(\Phi,Q)$, where $\Phi$ is an  $\BE$-stratifying system and $Q = \{Q_{i}\}_{i=1}^{n}$ is a set of indecomposables in $\A$ 
that are $\restr{\BE}{\F(\Phi)}$-projective in $\F(\Phi)$, 
such that for each $i$ there is an extriangle 
$
\begin{tikzcd}[column sep=0.5cm,cramped]
K_{i} \arrow{r}{} & Q_{i} \arrow{r}{q_{i}} & \Phi_{i} \arrow[dashed]{r}{\eta_{i}}& {},
\end{tikzcd}
$
satisfying some axioms (see \cref{def:weak-projective-system} for more details). 
The morphism $q_{i}$ then behaves like a projective cover of $\Phi_{i}$ in $\F(\Phi)$. 
If each $q_{i}$ is right minimal, then we call $(\Phi,Q)$ a \emph{minimal} projective $\BE$-stratifying system. 
We show the following.

\begin{thmx} 
\label{thmx:existence-of-proj-ESS}
\emph{(\cref{thm:existence-of-proj-ESS})}
Let  $\Phi$ be an $\BE$-stratifying system in an artin extriangulated category $(\A,\BE,\fs)$.
Then there is a set $Q$ of objects in $\F(\Phi)$, such that
$(\Phi,Q)$ is a minimal $\BE$-projective system. 
\end{thmx}

Another core notion we consider in this article is the \emph{Jordan-H{\"{o}}lder property} for filtrations; see \cref{def:E-s-notions}. 
Suppose $\Lambda$ is an arbitrary algebra, and that $\Phi$ is part of a projective stratifying system in the category of finite length left modules over $\Lambda$. 
Using the dual of \cite[Lem.\ 1.4]{ES}, it can be seen that the subcategory $\F(\Phi)$ 
has the Jordan-H{\"{o}}lder property on $\Phi$-filtrations. 
That is, any two $\Phi$-filtrations for an object $M\in\F(\Phi)$ are \emph{equivalent}, i.e.\ have the same length and the same ``composition factors'', up to permutation and isomorphism. 
Analogous phenomena occur for projective stratifying systems in the triangulated \cite[Prop.\ 5.11]{MS} and exact \cite[Prop.\ 6.5]{Sa} settings. 
The Jordan-H{\"{o}}lder property has been well studied for abelian categories, because the idea of \emph{length} behaves nicely in this context. Indeed, if an object in an abelian category has a composition series, then all such composition series are equivalent and, in particular, have the same length. 
For a general exact category $(\A,\E)$, the idea of length relative to the exact structure $\E$ has been studied in \cite{E19,BHLR}. 
It seems, however, that this kind of length is much better behaved when $(\A,\E)$ is \emph{Jordan-H{\"{o}}lder}, where one imposes that any two $\E$-composition series of an object must be equivalent (see \cite[Def.\ 5.1]{BHT}); compare, for example, \cite[Thm.\ 6.6]{BHLR} and \cite[Cor.\ 7.2]{BHT}.

Given a projective $\BE$-stratifying system $(\Phi,Q)$ in $(\A,\BE,\fs)$, we would like to show that $\F(\Phi)$ is (in some sense) a Jordan-H{\"{o}}lder extriangulated category. 
Such a notion would apply to triangulated categories as well. 
Unfortunately, the idea of length, and hence composition series, in a triangulated category is problematic due to 
the existence of triangles of the form 
$
\begin{tikzcd}[column sep=0.4cm,cramped]
X \arrow{r}& 0 \arrow{r}& X[1]\arrow{r}& X[1]
\end{tikzcd}
$ (see the discussion in \cref{subsec:grothendieck-monoid} following \cref{def:simplie-like-atom-like}).
This could be a major reason why the Jordan-H{\"{o}}lder property has only very recently been discussed in the general context of extriangulated categories; see Wang--Wei--Zhang--Zhang \cite{WWZZ}, also \cref{rem:WWZZ-simple}. 
In order to relate the Jordan-H{\"{o}}lder property on $\Phi$-filtrations in $\F(\Phi)$ to the functor $\restr{\BE}{\F(\Phi)}$, 
we introduce and discuss $(\BE,\fs)$-subobjects, 
-simple objects and 
-composition series; see \cref{sec:Grothendieck-monoid-and-JH-property}.  
We call an extriangulated category $(\A,\BE,\fs)$ \emph{Jordan-H{\"{o}}lder} if it has the Jordan-H{\"{o}}lder property on $(\BE,\fs)$-composition series, and \emph{length} if each object admits admits an $(\BE,\fs)$-composition series; see \cref{def:E-s-notions}. 
Our second main result is the following. 

\begin{thmx} 
\emph{(\cref{cor:JHP-for-F-Theta})}
Suppose $(\Phi,Q)$ is a minimal $\BE$-projective system in an artin extriangulated category  $(\A,\BE,\fs)$, such that  
$\A(Q_{i},-)$ is left exact on $(\F(\Phi),\restr{\BE}{\F(\Phi)},\restr{\fs}{\F(\Phi)})$ for each $Q_{i}\in Q$. 
Then 
$(\F(\Phi), \restr{\BE}{\F(\Phi)},\restr{\fs}{\F(\Phi)})$ 
is length and Jordan-H{\"{o}}lder. 
\end{thmx}

The final key notions we use are the Grothendieck group and monoid. 
Suppose $(\A,\E)$ is a skeletally small, exact category. 
The Jordan-H{\"{o}}lder property is closely linked to the \emph{Grothendieck group} $K_{0}(\A,\E)$ of $(\A,\E)$: 
it is straightforward to check that if $(\A,\E)$ is a Jordan-H{\"{o}}lder, length exact category, then $K_{0}(\A,\E)$ is a free abelian group. 
However, the converse is not true; see \cite[Exam.\ B]{E19}. 
This prompted Enomoto to study a related categorical invariant: 
the \emph{Grothendieck monoid} $M(\A,\E)$ of $(\A,\E)$ as introduced in \cite[Sec.\ 2.3]{BeGr}. 
It was shown that $(\A,\E)$ is a Jordan-H{\"{o}}lder, length exact category if and only if $M(\A,\E)$ is a free monoid (see \cite[Thm.\ 4.12]{E19}). 
Following Berenstein--Greenstein \cite{BeGr}, we extend the notion of a Grothendieck monoid to the extriangulated setting (see \cref{def:grothendieck-monoid}); see also independent work of Enomoto--Saito \cite{EnomotoSaito}. 
We then characterise the Jordan-H\"older property for extriangulated categories in terms of the Grothendieck monoid, which is our final main result.
We refer the reader to \cref{sec:Grothendieck-monoid-and-JH-property} for unexplained terms. 

\begin{thmx} 
\emph{(\cref{thm:JH-iff-free-monoid-iff-free-group})} 
The following are equivalent for a skeletally small, weakly idempotent complete extriangulated category $(\A,\BE,\fs)$. 
\begin{enumerate}[label=\textup{(\roman*)}]
    \item The extriangulated category $(\A,\BE,\fs)$ is length and Jordan-H{\"{o}}lder. 
    \item The Grothendieck monoid of $(\A,\BE,\fs)$ is free with basis in bijection with isoclasses of $(\BE, \fs)$-simple objects.
    \item The Grothendieck group of $(\A,\BE,\fs)$ is free with basis in bijection with isoclasses of $(\BE, \fs)$-simple objects.
\end{enumerate}
\end{thmx}

In \cref{sec:examples} we explore some examples. In particular, there are examples of projective stratifying systems in extriangulated categories that are not projective stratifying systems in exact or triangulated categories. We also give examples of Jordan-H{\"{o}}lder extriangulated categories: one arising from a stratifying system and one that does not.
Furthermore, Enomoto--Saito \cite{EnomotoSaito} recently asked if a skeletally small extriangulated category is exact if and only if it has a reduced Grothendieck monoid, and we give an explicit counterexample; see \cref{example:strong-proj-ESS} and \cref{rem:counter-to-EnomotoSaito-Q6-3}.

%%%%%%%%%%%%%%%%%%%%%%%%%%%%%%%%%%%%%%%%%%%%%%%%%
%%%%%%%%%%%%%%%%%%%%%%%%%%%%%%%%%%%%%%%%%%%%%%%%%
\section{Preliminaries}
\label{sec:background}

\subsection{Extriangulated categories}
\label{sec:extriangulated-categories}

Extriangulated categories were introduced in \cite{NP19} as simultaneous generalisations of exact and triangulated categories. In this section, we recall some details of these categories and useful terminology, but for further details we refer the reader to \cite[Sec.\ 2]{NP19}.
%
%%%%%%%%%%%%%%%%%%%%%%%%%%%%%%%%%%%%%%%%%%%%%%%%%
%
We denote the category of all abelian groups by $\Ab$. 

\begin{definition}
\label{def:extriangulated-category}
An \emph{extriangulated category} is a triplet 
$(\A, \mathbb{E}, \fs)$, where 
$\A$ is an additive category,    
$\mathbb{E}\colon \A^{\op} \times \A \to \Ab$ is a biadditive functor, and     
$\fs$ is an additive realisation of $\mathbb{E}$ (see \cite[Def.\ 2.10]{NP19}),
such that axioms (ET$3)$, (ET$3)^{\op}$, (ET$4)$ and (ET$4)^{\op}$ as stated in \cite[Def.\ 2.12]{NP19} are all satisfied. 
\end{definition}

Suppose $(\A,\BE,\fs)$ is an extriangulated category. For all pairs of objects $A,C\in\A$, an element $\delta\in\mathbb{E}(C,A)$ is called an \emph{extension}. 
The realisation $\fs$ (see \cite[Def.\ 2.9]{NP19}) assigns to each extension $\delta\in\mathbb{E}(C,A)$ an equivalence class 
\begin{equation}
\label{eqn:delta-realised-by-ABC}
\fs(\delta) = 
[\begin{tikzcd}
A \arrow{r}{a}& B\arrow{r}{b} & C 
\end{tikzcd}],
\end{equation}
where 
\begin{enumerate}[\textup{(\roman*)}]
    \item $\begin{tikzcd}[column sep=0.5cm,cramped]
    A \arrow{r}{a}& B\arrow{r}{b} & C 
    \end{tikzcd}$ 
    is a sequence of composable morphisms in $\A$; and 

    \item\label{equivalence} sequences 
    $\begin{tikzcd}[column sep=0.5cm,cramped]
    A \arrow{r}{a}& B\arrow{r}{b} & C 
    \end{tikzcd}$ 
    and 
    $\begin{tikzcd}[column sep=0.5cm,cramped]
    A \arrow{r}{a'}& B'\arrow{r}{b'} & C 
    \end{tikzcd}$ 
    in $\A$ are said to be \emph{equivalent} if there is an isomorphism 
    $g\colon B\to B'$ such that the diagram 
\begin{equation}
\label{eqn:equivalent-short-exact-sequences}
\begin{tikzcd}[row sep=0.6cm]
A \arrow{r}{a} \arrow[equals]{d}& B\arrow{r}{b} \arrow{d}{g}[swap]{\iso}& C \arrow[equals]{d}\\
A \arrow{r}{a'} & B'\arrow{r}{b'} & C 
\end{tikzcd}
\end{equation}
commutes (see \cite[Def.\ 2.7]{NP19}). 
    
\end{enumerate}

Let $\delta\in\mathbb{E}(C,A)$ be an extension. 
If \eqref{eqn:delta-realised-by-ABC} holds, then 
$a$ (resp.\ $b$) is called an \emph{inflation} (resp.\ \emph{deflation}), and 
the pair 
$
\langle
\begin{tikzcd}[column sep=0.5cm]
A \arrow{r}{a}& B\arrow{r}{b} & C, 
\end{tikzcd} 
\delta
\rangle
$ 
is called an \emph{extriangle}. 
Moreover, this is denoted by the diagram 
$
\begin{tikzcd}[column sep=0.5cm,cramped]
A \arrow{r}{a}& B\arrow{r}{b} & C \arrow[dashed]{r}{\delta}& {}. 
\end{tikzcd}
$
Given morphisms $x\colon A\to X$ and $z\colon Z\to C$, we obtain new extensions 
$x_{*}\delta \deff \BE(C,x)(\delta)$
and 
$x^{*}\delta \deff \BE(z,A)(\delta)$.
If $\eps\in\BE(C',A')$ is an extension, then a \emph{morphism of extensions} is a pair $(f,h)$, where $f\colon A\to A'$ and $h\colon C\to C'$ are morphisms in $\A$, such that $f_{*}\delta = h^{*}\eps$. 
By a \emph{morphism of extriangles}, we mean a commutative diagram
\[
\begin{tikzcd}[row sep=0.6cm]
A \arrow{r}{a} \arrow{d}{f}& B\arrow{r}{b} \arrow{d}{g}& C \arrow{d}{h} \arrow[dashed]{r}{\delta}& {}\\
A' \arrow{r}{a'} & B'\arrow{r}{b'} & C' \arrow[dashed]{r}{\eps}& {}
\end{tikzcd}
\]
where $(f,h)\colon \delta \to \eps$ is a morphism of extensions.

It follows from Yoneda's Lemma that there is a natural transformation
$\delta^{\sharp}\colon \A(A,-) \Rightarrow \BE(C,-)$
given by 
$(\delta^{\sharp})_{X}(x) 
    \deff x_{*}\delta
    = \mathbb{E}(C,x)(\delta)\in\mathbb{E}(C,X)$ 
for a morphism $x\colon A\to X$. 
There is a dually defined natural transformation 
$
\delta_{\sharp}\colon \A(-,C) \Rightarrow \BE(-,A)
$
(see \cite[Def.\ 3.1]{NP19}). 
More importantly, these fit into the exact sequences
\[
\hspace*{-0.1cm}
\begin{tikzcd}[column sep=1.07cm]
\A(C,Z) 
    \arrow{r}{\A(b,Z)} 
&\A(B,Z) 
    \arrow{r}{\A(a,Z)} 
&\A(A,Z) 
    \arrow{r}{(\delta^{\sharp})_{Z}} 
& \BE(C,Z) 
    \arrow{r}{\BE(b,Z)} 
&\BE(B,Z) 
    \arrow{r}{\BE(a,Z)} 
&\BE(A,Z)
\end{tikzcd}
\]
and 
\[ 
\hspace*{-0.1cm}
\begin{tikzcd}[column sep=1.07cm]
\A(Z,A) \arrow{r}{\A(Z,a)} &\A(Z,B) \arrow{r}{\A(Z,b)} &\A(Z,C) \arrow{r}{(\delta_{\sharp})_{Z}} 
& \BE(Z,A) \arrow{r}{\BE(Z,a)} &\BE(Z,B) \arrow{r}{\BE(Z,b)} &\BE(Z,C)
\end{tikzcd}
\]
in $\Ab$ for each $Z\in\A$ (see \cite[Cor.\ 3.12]{NP19}).

We recall some typical examples of extriangulated categories.

\begin{example}
\label{example:triangulated-is-extriangulated}
\cite[Prop.\ 3.22]{NP19} 
Suppose $\A$ is a triangulated category with suspension functor $[1]$. 
Then the triplet $(\A,\BE,\fs)$ is an extriangulated category, 
where: 
\begin{enumerate}[\textup{(\roman*)}]
    \item $\BE(-,-) = \A(-,-[1])$; and
    
    \item for 
    $\delta\in\BE(C,A) 
        = \A(C,A[1])$, 
    the realisation $\fs$ is defined as follows: 
    we complete $\delta$ to a triangle 
    $\begin{tikzcd}[column sep=0.5cm,cramped]
    A \arrow{r}{a}& B\arrow{r}{b} & C \arrow{r}{\delta}& A[1]  
    \end{tikzcd}$, then set 
    $\fs(\delta) = 
    [\begin{tikzcd}[column sep=0.5cm]
    A \arrow{r}{a}& B\arrow{r}{b} & C 
    \end{tikzcd}]$.
\end{enumerate} 
\end{example}

\begin{example}
\label{example:exact-is-extriangulated}
\cite[Exam.\ 2.13]{NP19} 
Suppose $\A$, equipped with a class $\E$ of short exact sequences, is an exact category. Let $A,C\in\A$, and set $\BE(C,A) \deff \Ext_{\E}(C,A)$ to be the collection of all sequences in $\E$ of the form 
$\begin{tikzcd}[column sep=0.5cm,cramped]
A \arrow{r}{a}& B\arrow{r}{b} & C 
\end{tikzcd}$ 
modulo the equivalence relation defined in \cref{def:extriangulated-category}\ref{equivalence}. 
In order to define a functor $\BE = \Ext_{\E} \colon \A^{\op}\times \A \to \Ab$, we must assume here that $\BE(C,A)$ is a set for all $A,C$. 
This is the case, for example, when $\A$ is skeletally small, or when $\A$ has enough projectives or enough injectives. 
With this set-theoretic assumption, we can use the existence of pushouts of inflations, pullbacks of deflations, and the Baer sum to ensure $\BE$ is a biadditive functor. The realisation $\fs$ is canonical: if 
$
\delta 
    = [\begin{tikzcd}[column sep=0.5cm]
    A \arrow{r}{a}& B\arrow{r}{b} & C 
    \end{tikzcd}]
$ 
is an extension set $\fs(\delta) = \delta$. 
It follows that 
$(\A,\BE,\fs)$ is an extriangulated category. 
\end{example}

We recall some notation from \cite[Def.\ 1.17]{LiuYNakaoka-hearts-of-twin-cotorsion-pairs-on-extriangulated-categories}. 
For any two classes $\X$, $\Y$ of objects of $\A$, we denote by $\X \ast \Y$ the class of all objects $A$ in $\A$ for which there exists an extriangle 
\[
\begin{tikzcd}
X \arrow{r}& A \arrow{r}& Y \arrow[dashed]{r}& {},
\end{tikzcd}
\]
with $X \in \X$ and $Y \in \Y$.
We consider $\X \ast \Y$ as a full subcategory of $\A$. 
Note also that the $\ast$ operation is associate: for classes $\X$, $\Y$ and $\Z$ of objects of $\A$, we have $(\X \ast \Y) \ast \Z = \X \ast (\Y \ast \Z)$; see \cite[Lem.\ 3.9]{Zh19}. 
Lastly, a full additive subcategory $\X\sse\A$ that is closed under isomorphisms is called \emph{extension-closed} if $\X \ast \X \sse \X$.

\begin{example}
\label{example:extension-closed-subcat-of-extriangulated-is-extriangulated}
It was shown in \cite{NP19} that an extension-closed subcategory $\B$ of an extriangulated category $(\A,\BE,\fs)$ is again extriangulated; see \cite[Def.\ 2.17, Rem.\ 2.18]{NP19}. 
Define $\restr{\BE}{\B}$ to be the restriction of $\BE$ to $\B^{\op}\times \B$. Since $\B$ is extension-closed, then we can restrict the realisation $\fs$ to $\restr{\BE}{\B}$ in order to obtain an additive realisation $\restr{\fs}{\B}$ of $\restr{\BE}{\B}$. 
In this way, $(\B, \restr{\BE}{\B},\restr{\fs}{\B})$ is an extriangulated category. 

A prominent example of this latter kind is discussed in \cite[Sec.\ 4.6]{PPPP}: Let $\Lambda$ be an Artin algebra. 
The homotopy category $K^b(\proj \Lambda)$ of bounded complexes of finitely generated projective modules over $\Lambda$ is a triangulated category, and hence extriangulated. 
The full subcategory $\A = K^{[-1,0]}(\proj \Lambda)$ of complexes concentrated in degrees $-1$ and $0$ is extension-closed, thus an extriangulated category. 
\end{example}

In the sequel, we will often concern ourselves with extriangulated categories $(\A, \mathbb{E}, \fs)$ whose underlying additive category $\A$ is \emph{Krull-Schmidt}, that is, each object in $\A$ decomposes into a finite direct sum of indecomposable objects
having local endomorphism rings and that this decomposition is unique up to isomorphism and permutation of summands. 
It follows from \cite[Cor.\ 4.4]{Krause2015} that every Krull-Schmidt category is \emph{weakly idempotent complete}, that is, every retraction admits a kernel. One may think of a weakly idempotent complete category as one that contains complements of direct summands that already exist in the category; see \cite[Appx.~A]{HenrardvanRoosmalen-Derived-categories-of-one-sided-exact-categories-and-their-localizations} or \cite[Sec.~2]{KlapprothMsapatoShah-idempotent-completion-of-n-exangulated-categories}, 
for more discussion. 
Note that any triangulated category is weakly idempotent complete.

We close this subsection with several results that will be useful in later sections. 
The first is a characterisation of weakly idempotent completeness through the lens of extriangulated category theory. 
For this, we recall a certain condition defined in \cite[Cond.~5.8]{NP19}. An extriangulated category $(\A,\BE,\fs)$ is said to satisfy the \emph{(WIC)} condition if, for composable morphisms $f,g$, we have $gf$ is an inflation (resp.\ deflation) implies $f$ is an inflation (resp.\ $g$ is a deflation).
The characterisation combines a result that appears in the PhD thesis \cite{Tattar-phd-thesis} of the fourth author and a result of Klapproth \cite{Klapproth-n-extension-closed-subcategories}.

\begin{prop}
\label{prop:weaklyidempotent}
Let $(\A, \BE, \fs)$ be an extriangulated category. Then the following statements are equivalent. 
\begin{enumerate}[label=\textup{(\roman*)}]
\item $\A$ is weakly idempotent complete.
\item Every retraction in $\A$ is a deflation in $(\A, \BE, \fs)$.
\item Every section in $\A$ is an inflation in $(\A, \BE, \fs)$. 
\item $(\A, \BE, \fs)$ satisfies (WIC).
\end{enumerate}
\end{prop}

\begin{proof} 
(i) $\Rightarrow$ (ii)\;\; 
Let $r\colon B \to C$ be a retraction with corresponding section $s\colon C \to B$ (so that $rs= \id{C}$) and let $k\colon K \to B$ be a kernel of $r$, which exists by assumption. 
It follows that $\begingroup
        \renewcommand\thickspace{\kern5pt}\bsm k & s \esm\endgroup \colon K\oplus C \to B$ is an isomorphism (see \cite[Rem.\ 7.4]{Buhler-exact-categories}). 
Then we have a commutative diagram
\[ 
\begin{tikzcd}[ampersand replacement=\&, column sep=1.3cm,row sep=0.6cm]
K 
    \arrow{r}{\begin{psmallmatrix}
    \id{K} \\ 0
    \end{psmallmatrix}}
    \arrow[d, equal] 
\& K \oplus C 
    \arrow{r}{
    \begingroup
        \renewcommand\thickspace{\kern5pt}
        \begin{psmallmatrix}
            0 & \id{C}
        \end{psmallmatrix}
    \endgroup}
    \arrow{d}{\iso}[swap]{
    \begingroup
        \renewcommand\thickspace{\kern5pt}
        \begin{psmallmatrix}
        k & s
        \end{psmallmatrix}
    \endgroup}
\& C 
    \arrow[d, equal] 
\\ 
K 
    \arrow[r, "k"'] 
\& B 
    \arrow[r, "r"'] 
\& C, 
\end{tikzcd} 
\] 
which is an equivalence of composable morphisms (see \eqref{eqn:equivalent-short-exact-sequences}). 
Our claim follows since the upper sequence is always an extriangle 
\cite[Rem.\ 2.11]{NP19}.
 
(ii) $\Rightarrow$ (i)\;\; 
Let $r\colon  B \to C$ be a retraction and suppose that there is an extriangle 
$
\begin{tikzcd}[column sep=0.4cm,cramped] 
A \arrow[r, "x"] & B \arrow[r, "r"] & C \arrow[r, dashed] & {} 
\end{tikzcd}
$ in $(\A, \BE, \fs)$. 
Note that this already implies $x$ is a weak kernel of $r$ by \cite[Prop.\ 3.3(2)]{NP19}. 
Let $s\colon C \to B$ be a morphism in $\A$ with $rs=\id{C}$. Then we apply (ET3) to obtain the morphism
\[ 
\begin{tikzcd}[ampersand replacement=\&, column sep=1.3cm,row sep=0.6cm]
A 
    \arrow[r, "\bsm \id{A} \\ 0 \esm"] 
    \arrow[d, equal] 
\& A \oplus C 
    \arrow{r}{
    \begingroup
        \renewcommand\thickspace{\kern5pt}
        \begin{psmallmatrix}
        0 & \id{C} 
        \end{psmallmatrix}
    \endgroup} 
    \arrow{d}{\begingroup
        \renewcommand\thickspace{\kern5pt}
        \begin{psmallmatrix} 
            x & s 
        \end{psmallmatrix}
    \endgroup}
\& C 
    \arrow[dotted]{d}{\exists c}
    \arrow[r, dashed] 
\& {}  
\\ 
A 
    \arrow[r, "x"'] 
\& B 
    \arrow[r, "r"'] 
\& C 
    \arrow[r, dashed] 
\& {}
\end{tikzcd} 
\] 
of extriangles. By the commutativity of the right hand square, we have that 
$
\begingroup
        \renewcommand\thickspace{\kern5pt}
        \bsm 0 & c \esm\endgroup 
    = r \begingroup
        \renewcommand\thickspace{\kern5pt}
        \bsm x & s \esm\endgroup  
    = \begingroup
        \renewcommand\thickspace{\kern5pt}
        \bsm 0 & \id{C} \esm\endgroup
$, whence $c = \id{C}$. 
It follows from \cite[Cor.\ 3.6]{NP19} that 
$\begingroup
        \renewcommand\thickspace{\kern5pt}
        \bsm x & s \esm\endgroup
        \colon  A \oplus C \to B
$ 
is an isomorphism, and we deduce that $x$ is a monomorphism as the composite of two monomorphisms. Thus, $x$ is a kernel of $r$. 

The equivalence of (i) and (iii) is similar. The equivalence of (i) and (iv) is \cite[Prop.~2.7]{Klapproth-n-extension-closed-subcategories}.
\end{proof}

The following definition is a special case of \cite[Def.\ 4.31]{KlapprothMsapatoShah-idempotent-completion-of-n-exangulated-categories}.

\begin{definition}
\label{def:WIC-extriangulated-category}
An extriangulated category $(\A, \BE, \fs)$ is called \emph{weakly idempotent complete} if the underlying additive category $\A$ is weakly idempotent complete in the usual sense, and hence $(\A, \BE, \fs)$ satisfies the equivalent conditions from \cref{prop:weaklyidempotent}.
\end{definition}

For the remaining technical results, we recall the additive structure on the category of extriangles.
Let $A,A',C,C'$ be objects of an extriangulated category $(\A,\BE,\fs)$. We denote by ${}_{A}0_{C}$ the split extension (i.e.\ the abelian group identity) in $\BE(C,A)$. 
Let $\delta\in\BE(C,A)$ and $\delta'\in\BE(C',A')$ be any extensions. 
Following \cite[Def.\ 2.6]{NP19}, we use the notation $\delta\oplus\delta'$ for the extension in $\BE(C\oplus C',A\oplus A')$ corresponding to the element $(\delta,0,0,\delta')$ under the natural isomorphism
\begin{equation}
\label{eqn:bifunctor-isomorphism}
\BE(C\oplus C',A\oplus A') \iso 
    \BE(C,A) \oplus 
    \BE(C,A') \oplus 
    \BE(C',A) \oplus 
    \BE(C',A'). 
\end{equation}
The \emph{category of $\BE$-extensions} $\BE\dExt{\A}$ has as its objects all the $\BE$-extensions, and as its morphisms the morphisms of extensions in the usual sense; see \cite[Sec.\ 2.1]{NP19}.
It can be shown that $\BE\dExt{\A}$ is an additive category with this direct sum of extensions; 
see 
\cite[Prop.\ 3.2]{Bennett-TennenhausHauglandSandoyShah-the-category-of-extensions-and-a-characterisation-of-n-exangulated-functors}. 
Furthermore, this operation induces an additive structure on the category of extriangles (see \cite[Rem.\ 2.14, Prop.\ 3.3]{HerschendLiuNakaoka-n-exangulated-categories-I-definitions-and-fundamental-properties}) as follows. 
If
$
\langle
\begin{tikzcd}[column sep=0.5cm]
A \arrow{r}{a}& B\arrow{r}{b} & C, 
\end{tikzcd} 
\delta
\rangle
$ 
and 
$
\langle
\begin{tikzcd}[column sep=0.5cm]
A' \arrow{r}{a'}& B'\arrow{r}{b'} & C', 
\end{tikzcd} 
\delta'
\rangle
$
are extriangles, then their direct sum is given by 
$
\langle
\begin{tikzcd}%[column sep=0.5cm]
A\oplus A' \arrow{r}{a\oplus a'}& B\oplus B'\arrow{r}{b\oplus b'} & C\oplus C', 
\end{tikzcd} 
\delta \oplus \delta'
\rangle
$.

The following lemmas yield a generalised version of \cite[Lem.~2.2]{adachi2022mixed} to the weakly idempotent complete setting.

\begin{lem}
\label{lem:decomposable}
Let $(\A, \BE, \fs)$ be a weakly idempotent complete extriangulated category. 
Then any extriangle of the form
\[
\begin{tikzcd}[ampersand replacement=\&, column sep=1.2cm]
A 
    \arrow{r}{\begin{psmallmatrix}
       a_{1} \\ a_{2}
    \end{psmallmatrix}}
\& B_{1} \oplus B_{2} 
    \arrow{r}{\begingroup
        \renewcommand\thickspace{\kern5pt}
        \begin{psmallmatrix}
        0 &b_{2}
        \end{psmallmatrix}
        \endgroup}
\& C 
    \arrow[r, dashed, "\delta"] 
\&{}
\end{tikzcd} 
\] 
is isomorphic to 
\[
\langle
\begin{tikzcd}[column sep=1cm]
B_{1}\arrow{r}{\id{B_{1}}}& B_{1}\arrow{r}{} & 0, 
\end{tikzcd} 
{}_{B_{1}}0_{0}
\rangle
\oplus 
\langle
\begin{tikzcd}[column sep=1cm]
A' \arrow{r}{a_{2}s}& B_{2}\arrow{r}{b_{2}} & C, 
\end{tikzcd} 
\delta'
\rangle,
\]
where
there is an isomorphism 
$
\begingroup
        \renewcommand\thickspace{\kern5pt}
        \bsm x & s \esm
\endgroup
    \colon 
    B_{1} \oplus A' \overset{\iso}{\longrightarrow}
    A
$.
%In particular, $b_{2}$ is also a deflation.
\end{lem}
\begin{proof}
Since 
$\begingroup
\renewcommand\thickspace{\kern5pt}
\bsm 0 &b_{2} \esm
\endgroup
	= b_{2}\circ \begingroup
\renewcommand\thickspace{\kern5pt}
\bsm  0 & \id{B_{2}} \esm
\endgroup
$ is a deflation, 
we know that $b_{2}$ is a deflation by (WIC) as $(\A, \BE, \fs)$ is weakly idempotent complete (see \cref{prop:weaklyidempotent}). 
Thus, there is an extriangle of the form
$
\begin{tikzcd}[column sep=0.5cm,cramped]
A' \arrow{r}{a}& B_{2}\arrow{r}{b_{2}} & C \arrow[dashed]{r}{\delta'}&{.}
\end{tikzcd} 
$ 
Taking the direct sum with 
$
\begin{tikzcd}[column sep=0.8cm,cramped]
B_{1}\arrow{r}{\id{B_{1}}}& B_{1}\arrow{r}{} & 0 \arrow[dashed]{r}{{}_{B_{1}}0_{0}}&{}
\end{tikzcd}
$
and applying (ET$3)^{\op}$, there is a morphism of extriangles
\[
\begin{tikzcd}[ampersand replacement=\&, column sep=1.5cm]
B_{1} \oplus A' 
	\arrow{r}{\begingroup
        \renewcommand\thickspace{\kern5pt}
        \begin{psmallmatrix}
        \id{B_{1}} & 0 \\
        0 & a
        \end{psmallmatrix}
        \endgroup}
	\arrow[dotted]{d}{\iso}[swap]{\begingroup
        \renewcommand\thickspace{\kern5pt}
        \begin{psmallmatrix}
        x & s
        \end{psmallmatrix}
        \endgroup}
\& B_{1} \oplus B_{2}
	\arrow{r}{\begingroup
        \renewcommand\thickspace{\kern5pt}
        \begin{psmallmatrix}
        0 & 0 \\
        0 &b_{2}
        \end{psmallmatrix}
        \endgroup}
	\arrow[equals]{d}
\& 0\oplus C 
	\arrow{d}{\iso}[swap]{\begingroup
        \renewcommand\thickspace{\kern5pt}
        \begin{psmallmatrix}
        0 \\
        \id{C} 
        \end{psmallmatrix}
        \endgroup}
	\arrow[dashed]{r}{{}_{B_{1}}0_{0} \oplus \delta'}
\& {}
\\
A 
	\arrow{r}{\begingroup
        \renewcommand\thickspace{\kern5pt}
        \begin{psmallmatrix}
        a_{1} \\
        a_{2}
        \end{psmallmatrix}
        \endgroup}
\& B_{1} \oplus B_{2} 
	\arrow{r}{\begingroup
        \renewcommand\thickspace{\kern5pt}
        \begin{psmallmatrix}
        0 & b_{2}
        \end{psmallmatrix}
        \endgroup}
\& C
	\arrow[dashed]{r}{\delta}
\& {}
\end{tikzcd}
\]
where $(x,s)$ is an isomorphism by \cite[Cor.~3.6]{NP19}. 
By commutativity of the diagram above, we see $a = a_{2}s$ and we are done. 
\end{proof}

\begin{lem}\label{lem:decomposable3}
Let $(\A, \BE, \fs)$ be a weakly idempotent complete extriangulated category. 
Then any extriangle of the form
\[
\begin{tikzcd}[ampersand replacement=\&, column sep=1.5cm]
A 
    \arrow{r}
\& B_{1} \oplus B_{2} 
    \arrow{r}{\begingroup
        \renewcommand\thickspace{\kern5pt}
        \begin{psmallmatrix}
        g_{1} & 0\\
        0 & g_{2}
        \end{psmallmatrix}
        \endgroup}
\& C_{1} \oplus C_{2} 
    \arrow[r, dashed, "\delta"] 
\&{}
\end{tikzcd} 
\] 
is isomorphic to 
$
\langle
\begin{tikzcd}[ampersand replacement=\&]
 A_{1}
    \arrow{r}{f_{1}}
\& B_{1} 
    \arrow{r}{g_{1}}
    \& C_{1},
\end{tikzcd} 
\delta_{1}
\rangle
\oplus
\langle
\begin{tikzcd}[ampersand replacement=\&]
 A_{2}
    \arrow{r}{f_{2}}
\& B_{2} 
    \arrow{r}{g_{2}}
    \& C_{2},
\end{tikzcd} 
\delta_{2}
\rangle.
$
\end{lem}

\begin{proof}
Since $\fs$ is an additive realisation, we know that 
$ \begingroup
\renewcommand\thickspace{\kern5pt}
\bsm  0 & \id{C_{2}} \esm
\endgroup
\colon C_{1}\oplus C_{2} \to C_{2}$
is a deflation. 
Since 
$\begingroup
\renewcommand\thickspace{\kern5pt}
\bsm 0 & g_{2} \esm
\endgroup
	= \begingroup
\renewcommand\thickspace{\kern5pt}
\bsm  0 & \id{C_{2}} \esm
\endgroup
\circ 
\begingroup
        \renewcommand\thickspace{\kern5pt}
        \begin{psmallmatrix}
        g_{1} & 0\\
        0 & g_{2}
        \end{psmallmatrix}
        \endgroup
$
is the composition of deflations, it is a deflation itself by \cite[Rem.~2.16]{NP19}. 
Like in the proof of \cref{lem:decomposable}, it then follows that $g_{2}$ is a deflation and we have an extriangle of the form
$
\begin{tikzcd}[ampersand replacement=\&,cramped]
 A_{1}
    \arrow{r}{f_{1}}
\& B_{1} 
    \arrow{r}{g_{1}}
    \& C_{1}
    		\arrow[dashed]{r}{\delta_{1}}
    \&{.}
\end{tikzcd} 
$
By symmetry, we also have an extriangle 
$
\begin{tikzcd}[ampersand replacement=\&,cramped]
 A_{2}
    \arrow{r}{f_{2}}
\& B_{2} 
    \arrow{r}{g_{2}}
    \& C_{2}
    		\arrow[dashed]{r}{\delta_{2}}
    \&{.}
\end{tikzcd} 
$
Using (ET$3)^{\op}$ and \cite[Cor.~3.6]{NP19}, we obtain an isomorphism of extriangles 
\[
\begin{tikzcd}[ampersand replacement=\&, column sep=1.5cm]
A 
	\arrow[dotted]{d}{\iso}
    \arrow{r}
\& B_{1} \oplus B_{2} 
	\arrow[equals]{d}
    \arrow{r}{g_{1}\oplus g_{2}}
\& C_{1} \oplus C_{2} 
    \arrow[r, dashed, "\delta"] 
    	\arrow[equals]{d}
\&{} \\
A_{1} \oplus A_{2}
    \arrow{r}{f_{1} \oplus f_{2}}
\& B_{1} \oplus B_{2}
    \arrow{r}{g_{1} \oplus g_{2}}
    \& C_{1} \oplus C_{2}
    		\arrow[dashed]{r}{\delta_{1} \oplus \delta_{2}}
    \&{,}
\end{tikzcd} 
\]
which completes the proof.
\end{proof}

The last result of this subsection is an analogue of \cite[Prop.~2.1]{IyamaYoshino-mutation-in-tri-cats-rigid-CM-mods} for a Krull-Schmidt extriangulated category, which was first shown in \cite[Lem.~2.5(2)]{adachi2022mixed}. We include full details of the proof for convenience.

\begin{lem} \label{lem:summandslemma}
Let $(\A, \BE, \fs)$ be a Krull-Schmidt extriangulated category
Suppose $\X, \Y$ are subcategories of $\A$ that are closed under direct summands. If $\A(\X, \Y) = 0$, then $\X \ast \Y$ is closed under direct summands.
\end{lem}

\begin{proof}
Assume that $Z \deff Z_{1} \oplus Z_{2} \in \X \ast \Y$. 
Since $Z\in\X\ast\Y$, there is an extriangle
\begin{equation}\label{eqn:KS-extriangle-for-Z}
\begin{tikzcd}%[column sep=0.5cm] 
X \arrow{r}{f}& Z \arrow{r}{g}& Y \arrow[dashed]{r}{}& {,}
\end{tikzcd}
\end{equation}
where $X\in\X$ and $Y\in\Y$. Note that $f$ is a right $\X$-approximation of $Z$ since $\A(\X,\Y) = 0$. 

If $f=0$, then by the dual of \cref{lem:decomposable}, we have that \eqref{eqn:KS-extriangle-for-Z} is isomorphic to 
\[
\langle
\begin{tikzcd}[column sep=0.7cm]
0
	\arrow{r}
&Z_{1}
	\arrow{r}{\id{Z_{1}}}
& Z_{1},
\end{tikzcd} 
{}_{0}0_{B_{1}}
\rangle
\oplus 
\langle
\begin{tikzcd}[column sep=0.7cm]
X
	\arrow{r}{}
& Z_{2}
	\arrow{r}{} 
& Y', 
\end{tikzcd} 
\delta'
\rangle,
\]
where $Z_{1} \oplus Y' \iso Y$. 
We see that $Z_{1}\in\Y$ because $\Y$ is closed under direct summands. By symmetry, we also have $Z_{2}\in\Y$. Noting that $\Y\sse\X\ast\Y$ concludes the case when $f=0$. 

Thus, now assume $f\neq 0$. By \cite[Cor.~1.4]{KrauseSaorin-minimal-approximations-of-modules}, we may assume that $f$ is of the form
$\begingroup
        \renewcommand\thickspace{\kern5pt}
        \bsm f' & f'' \esm
\endgroup
        \colon
X'\oplus X'' \to Z$, where $f'\neq 0$ is right minimal and $f''=0$. 
So, as $f = \begingroup
        \renewcommand\thickspace{\kern5pt}
        \bsm f' & 0 \esm
\endgroup$
and 
$\begingroup
        \renewcommand\thickspace{\kern5pt}
        \bsm \id{X'} \\ 0 \esm
\endgroup\colon X' \to X'\oplus X''$
are inflations, 
we see that $f' = f\circ\begingroup
        \renewcommand\thickspace{\kern5pt}
        \bsm \id{X'} \\ 0 \esm
\endgroup $
is an inflation by \cite[Rem.~2.16]{NP19}. 
Furthermore, since $\X$ is closed under direct summands, we have $X'\in\X$ and $f'\colon X'\to Z$ is a minimal right $\X$-approximation of $Z$. 
Moreover, we have that there is an extriangle
\begin{equation}\label{eqn:minimal-extriangle-for-Z}
\begin{tikzcd}%[column sep=0.5cm] 
X' \arrow{r}{f'}& Z \arrow{r}{g'}& Y' \arrow[dashed]{r}{}& {,}
\end{tikzcd}
\end{equation}
where $Y'$ is a direct summand of $Y$ and hence $Y'\in\Y$. 
%%%% Proof that Y' is a summand of Y
%By (ET3) and (ET3$)^{\op}$ we have morphisms of extriagles
%\[
%\begin{tikzcd}
%X' \oplus X'' 
%	\arrow{r}{f}
%	\arrow{d}{(1,0)}
%& Z 
%	\arrow{r}{}
%	\arrow[equals]{d}
%& Y 
%	\arrow[dashed]{r}{}
%	\arrow[bend right]{d}{y}
%& {} \\
%X' 
%	\arrow{r}{f'}
%	\arrow{u}{(1,0)^T}
%& Z
%	\arrow{r}{f'}
%& Y'
%	\arrow[dashed]{r}{}
%	\arrow[bend left]{u}{r}
%& {}
%\end{tikzcd}
%\]
%Then $yr$ is an automorphism of $Y'$, so $y$ is a retraction and hence $Y'$ is a summand of $Y$. 
%%%%

Note that the compositions
$\begingroup
\renewcommand\thickspace{\kern5pt}
\bsm \id{Z_{1}} & 0 \esm
\endgroup\circ f' \colon X'\to Z_{1}$
and 
$\begingroup
\renewcommand\thickspace{\kern5pt}
\bsm 0& \id{Z_{2}} \esm
\endgroup\circ f' \colon X'\to Z_{2}$
are both right $\X$-approximations. 
First, assume 
$\begingroup
\renewcommand\thickspace{\kern5pt}
\bsm \id{Z_{1}} & 0 \esm
\endgroup f' = 0$.
Then we may apply the dual of \cref{lem:decomposable} to \eqref{eqn:minimal-extriangle-for-Z} to obtain an extriangle
$\begin{tikzcd}[column sep=0.5cm,cramped] 
X' \arrow{r}{}& Z_{2} \arrow{r}{}& Y'' \arrow[dashed]{r}{}& {,}
\end{tikzcd}$
where $Z_{1} \oplus Y'' \iso Y'$. In this case, we see that 
$Z_{1}\in\Y$ and 
$Z_{2}\in\X\ast\Y$, and we are done. 
By symmetry, we are also done if $\begingroup
\renewcommand\thickspace{\kern5pt}
\bsm 0& \id{Z_{2}} \esm
\endgroup f' = 0$.

Therefore, assume now that both 
$\begingroup
\renewcommand\thickspace{\kern5pt}
\bsm \id{Z_{1}} & 0 \esm
\endgroup f'$
 and 
 $\begingroup
\renewcommand\thickspace{\kern5pt}
\bsm 0& \id{Z_{2}} \esm
\endgroup f'$ 
are non-zero. 
In particular, we may take minimal right $\X$-approximations 
$f_{1}\colon X_{1} \to Z_{1}$
and 
$f_{2}\colon X_{2} \to Z_{2}$ of 
$Z_{1}$ and $Z_{2}$, respectively.
This yields a minimal right $\X$-approximation 
$f_{1}\oplus f_{2} \colon X_{1} \oplus X_{2} \to Z_{1}\oplus Z_{2}$
of $Z = Z_{1} \oplus Z_{2}$. 
(Indeed, since $\A$ is Krull-Schmidt, by \cite[Cor.~1.4]{KrauseSaorin-minimal-approximations-of-modules} there is a direct sum decomposition of $X_1 \oplus X_2 \cong A \oplus B$ such that $\restr{(f_1 \oplus f_2)}{A}$ is right minimal and $\restr{(f_1 \oplus f_2)}{B} = 0$, so if $B \neq 0$ then this would contradicts the minimality of $f_1$ or $f_2$.)
This implies  $f_{1}\oplus f_{2}$ is isomorphic to $f'$ by uniqueness of minimal approximations (see e.g.\ \cite[Prop.~2.2]{AuslanderReiten-Rep-theory-of-Artin-algebras-IV}), 
and hence
$f_{1}\oplus f_{2}$ is an inflation. 
Using \cref{lem:decomposable3}, we see that there are extriangles
\[
\begin{tikzcd}[ampersand replacement=\&]
 X_{1}
    \arrow{r}{f_{1}}
\& Z_{1} 
    \arrow{r}{g_{1}}
    \& Y_{1}
    \arrow[dashed]{r}{\delta_{1}}
    \& {}
\end{tikzcd} 
\hspace{0.5cm}\text{and}\hspace{0.5cm}
\begin{tikzcd}[ampersand replacement=\&]
 X_{2}
    \arrow{r}{f_{2}}
\& Z_{2} 
    \arrow{r}{g_{2}}
    \& Y_{2}
    	\arrow[dashed]{r}{\delta_{2}}
    \& {,}
\end{tikzcd} 
\]
where $Y_{1} \oplus Y_{2} \iso Y'$. As $\Y$ is closed under direct summands, this implies $Z_{i}\in\X\ast\Y$ for $i=1,2$.
\end{proof}

%%%%%%%%%%%%%%%%%%%%%%%%%%%%%%%%%%%%%%%%%%%%%%%%%

\subsection{Filtered objects}
\label{subsec:filtered-objects-in-an-extriangulated-category}

For this subsection, let $(\A,\BE,\fs)$ denote an extriangulated category and let $\X$ be a (possibly empty) class of objects of $\A$. 
We discuss the notion of $\X$-filtered objects 
in $(\A,\BE,\fs)$ as introduced in \cite{Zh19}, which generalises the corresponding notion for exact and triangulated categories found in \cite{Sa} and \cite{MS}, respectively.

\begin{definition}
\label{def:X-filtration}
Let $M\in\A$ be an object. 
An \emph{$\X$-filtration of $M$} 
is a 
sequence 
$\xi = (\xi_{0},\xi_{1},\ldots,\xi_{t})$ 
of extriangles 
$
\begin{tikzcd}[column sep=0.5cm,cramped]
M_{i-1} \arrow{r} & M_{i} \arrow{r} & X_{i}
\arrow[dashed]{r}{\xi_{i}} & {}
\end{tikzcd}
$
with 
$M_{-1} = M_{0} = X_{0} = 0$, 
$M_{t} = M$, 
and 
$X_{i} \in \X$ 
for all $1\leq i\leq t$. 
We call $t$ the \emph{length} of $\xi$ and 
such an $\X$-filtration of $M$ is depicted as follows.
\begin{equation}
\label{eqn:X-filtration}
\begin{tikzcd}[column sep=0.5cm, row sep=0.5cm]
0 
% = M_{-1} 
\arrow{rr}&&
0 
% = M_{0} 
\arrow{rr}\arrow{dl}&& M_{1} \arrow{rr}\arrow{dl}&& M_{2} \arrow{r}\arrow{dl}& \cdots \arrow{r}& M_{t-1} \arrow{rr}&& M_{t} = M \arrow{dl}\\
& 0
% = X_{1}
\arrow[dashed]{ul}{\xi_{0}}&&X_{1}\arrow[dashed]{ul}{\xi_{1}}&&X_{2}\arrow[dashed]{ul}{\xi_{2}} &&&&X_{t} \arrow[dashed]{ul}{\xi_{t}}&
\end{tikzcd}
\end{equation}

We say $M$ is \emph{$\X$-filtered} if 
it admits an $\X$-filtration. 
We denote by $\F(\X)$ the full subcategory of $\A$ formed by all $\X$-filtered objects of $\A$. 
\end{definition}

\begin{remark}
\label{rem:comments-on-FTheta}
Let us make some observations about \cref{def:X-filtration}.
\begin{enumerate}[\textup{(\roman*)}]
    
    \item Any object $M$ which has an $\X$-filtration of length $t= 0$ is a zero object.
    
    \item\label{part:F-empty-is-zero} If $\X = \emptyset$ is empty, then we have that $\F(\X) = \{0\}$.

    \item\label{part:F(X)-is-smallest-ext-closed-subcat-containing-X} 
    The subcategory $\F(\X)$ is the smallest extension-closed subcategory of $(\A, \BE,\fs)$ containing $\X$ \cite[Lem.\ 3.2]{Zh19}. 
    It follows from 
    \cref{example:extension-closed-subcat-of-extriangulated-is-extriangulated} that 
    $(\F(\X), \restr{\BE}{\F(\X)},\restr{\fs}{\F(\X)})$ 
    is an extriangulated category. 
    In fact, it is an extriangulated subcategory of $(\A, \BE,\fs)$ in the sense of Haugland \cite[Def.\ 3.7]{Haugland-the-grothendieck-group-of-an-n-exangulated-category} and canonical inclusion functor $\F(\X) \into \A$ forms part of an extriangulated functor as defined in \cite[Def.\ 2.32]{Bennett-TennenhausShah-transport-of-structure-in-higher-homological-algebra}.
    
    \item The objects $X_{1},\ldots,X_{t}$ in \eqref{eqn:X-filtration} can be thought of as ``composition factors'' of $M$ with respect to this specific $\X$-filtration, and we refer to them as the \emph{$\X$-factors} or just the \emph{factors}. However, it is not clear that these are unique up to isomorphism and permutation, i.e.\ that $\F(\X)$ has the Jordan-H{\"{o}}lder property. 
    We show in Section~\ref{sec:stratifying-systems} that $\F(\X)$ satisfies the Jordan-H{\"{o}}lder property whenever $\X$ forms part of what we call a projective $\BE$-stratifying system satisfying some mild assumption (see \cref{thm:JHP-for-Theta-filtrations}).

\end{enumerate}
\end{remark}

Given subcategories $\Y$ and $\Z$ of $\A$, we use the notation $\A(\Y,\Z) = 0$ to mean $\A(Y,Z) = 0$ for all $Y\in\Y$ and $Z\in\Z$. 
Similar notation is defined analogously. 
We use the following extriangulated version of \cite[Lem.\ 4.3]{MS} several times.

\begin{lem}
\label{lem:Zhou-lem3-3-orthogonality-compatible-with-F}
\cite[Lem.\ 3.3]{Zh19} 
Let $\Y$ and $\Z$ be full subcategories of $(\A,\BE,\fs)$.  
\begin{enumerate}[\textup{(\roman*)}]
    \item\label{part:Zhou3.3i} If $\A(\Y,\Z) = 0$, then $\A(\F(\Y),\F(\Z)) = 0$.
    \item\label{part:Zhou3.3ii} If $\BE(\Y,\Z) = 0$, then $\BE(\F(\Y),\F(\Z)) = 0$.
\end{enumerate}
\end{lem}

% %%%%%%%%%%%%%%%%%%%%%%%%%%%%%%%%%%%%%%

\section{The Grothendieck monoid and the Jordan-H{\"{o}}lder property} 
\label{sec:Grothendieck-monoid-and-JH-property}

For an exact category, a characterisation of the Jordan-H{\"{o}}lder property in terms of Grothendieck monoids is given in \cite[Sec.\ 4.3]{E19}. 
Our goal in this section is to prove an analogue for extriangulated categories. 
In Subsection~\ref{subsec:JHP-defs} we motivate and define fundamental notions we will work with, and in Subsection~\ref{subsec:grothendieck-monoid} we introduce and study the Grothendieck monoid of an extriangulated category. This allows us to characterise length Jordan-H{\"{o}}lder extriangulated categories using the language of Grothendieck monoids in Subsection~\ref{subsec:characterisation-of-JHP}. 
As noted in the introduction, some notions have been introduced independently in \cite{EnomotoSaito,WWZZ}.

%%%%%%%%%%%%%%%%%%%%%%%%%%%%%%%%%%%%%%%%%%%%%%%%%%
%%%%%%%%%%%%%%%%%%%%%%%%%%%%%%%%%%%%%%%%%%%%%%%%%%

\subsection{Jordan-H{\"{o}}lder extriangulated categories}
\label{subsec:JHP-defs}

To discuss some version of the Jordan-H{\"{o}}lder property for extriangulated categories, we need a notion of length. 
Following ideas from the triangulated and exact cases, Zhou \cite[Sec.\ 3]{Zh19} defines a length on an object $M\in\F(\X)$ with respect to the class of objects $\X$ in an extriangulated category $(\A,\BE,\fs)$; 
see also \cite[Sec.\ 4]{MS}, \cite[Sec.\ 5]{Sa}. 
The \emph{$\X$-length $\ell_{\X}(M)$} of $M$ is defined to be the minimum of the set 
\[
\Set{
t\in\bN 
|
\exists \X\text{-filtration } (\xi_{1},\ldots,\xi_{t}) \text{ of } M 
}.
\]
We note that there are parallels between the definition of $\X$-length 
and the dimension of a triangulated category as defined in \cite[Def.\ 3.2]{Rou}.

There is an immediate limitation of the $\X$-length: it depends entirely on the class of objects $\X\sse\A$. 
Indeed, if we take $\X = \A$, then every non-zero object has $\X$-length one---even 
$\ell_{\X}(M) = \ell_{\X}(M\oplus M)$---which is not desirable. 
This illustrates the need to define a length 
which depends only on the extriangulated structure.  In this paper,  we will show that if the class $\X$ satisfies some conditions, then one may define a sensible $\X$-length on the subcategory of $\X$-filtered objects.

Throughout this subsection, let $(\A,\BE,\fs)$ denote an extriangulated category.
Part of the following definition generalises \cite[Defs.\ 3.1, 3.3]{BHLR}.

\begin{definition}
\label{def:E-s-notions}
\begin{enumerate}[label=\textup{(\roman*)}]
    \item\label{item:E-s-subobject} For objects $A,B\in\A$, 
    we say that $A$ is an
    \emph{$(\BE,\fs)$-subobject of $B$} 
    if there is an extriangle of the form 
    $
    \begin{tikzcd}[column sep=0.5cm,cramped]
    A\arrow{r}{a}&B\arrow{r}&C\arrow[dashed]{r}{}&{}
    \end{tikzcd}
    $
    in $(\A,\BE,\fs)$. 
    
    \item\label{item:E-s-simple} 
    A non-zero object $S\in\A$ is called \emph{$(\BE,\fs)$-simple} if, whenever there is an extriangle 
    $
    \begin{tikzcd}[column sep=0.5cm,cramped]
        A\arrow{r}{a}
        &S\arrow{r}
        &C\arrow[dashed]{r}{}
        &{},
    \end{tikzcd}
    $ 
    either 
    $a$ = 0 or $a$ is an isomorphism.

    \item\label{item:E-s-composition-series}
    We denote by 
    $\Sim(\A,\BE,\fs)$ 
    the class of all $(\BE,\fs)$-simple objects of $(\A,\BE,\fs)$. 
    We call an $\Sim(\A,\BE,\fs)$-filtration (in the sense of \cref{def:X-filtration}) of an object $M \in (\A,\BE,\fs)$ 
    an \emph{$(\BE,\fs)$-composition series}. Note these are finite by definition. 
    
    \item\label{item:E-s-length}
    We say that $(\A,\BE,\fs)$ is a 
    \emph{length} extriangulated category 
    if each 
    object 
    $M\in\A$ has an $(\BE,\fs)$-composition series 
    and the set of lengths of all $(\BE,\fs)$-composition series of $M$ is bounded.
    
    \item \label{item:Jordan-Holder-extrianguated-category}
We call two $(\BE,\fs)$-composition series 
$(\xi_{0},\ldots,\xi_{t})$ 
and 
$(\xi'_{0},\ldots,\xi'_{s})$ 
of an object 
\emph{equivalent} 
if 
$s=t$ and, up to a permutation of indices, 
the $(\BE,\fs)$-simple factors are isomorphic. 
We say $(\A,\BE,\fs)$ is a \emph{Jordan-H{\"{o}}lder} extriangulated category 
if, for each object $M\in\A$, any two  $(\BE,\fs)$-composition series of $M$ are equivalent. 
\end{enumerate}
\end{definition}

\begin{remark}
\label{rem:WWZZ-simple}
Our definition of a simple object with respect to the extriangulated structure $(\BE,\fs)$ is more general than what is called `simple' in \cite[Sec.\ 3]{WWZZ}.
Indeed, our version allows for extriangles of the form 
$
    \begin{tikzcd}[column sep=0.5cm,cramped]
        A\arrow{r}{0}
        &S\arrow{r}
        &C\arrow[dashed]{r}{}
        &{},
    \end{tikzcd}
$ 
which are commonplace in triangulated categories. 
However, both versions recover the notion of a simple object in an abelian category or in an exact category.

Accordingly, an $(\BE,\fs)$-composition series (in the sense of \cref{def:E-s-notions}) is more general than a `composition series' in the sense of Wang--Wei--Zhang--Zhang. An extriangulated category is Jordan-H{\"{o}}lder and length in our sense if and only if it satisfies the `Jordan-H{\"{o}}lder property' defined in \cite[Def.\ 3.1]{WWZZ} if one chooses $\X$ to be the collection of $(\BE,\fs)$-simple objects. 
\end{remark}

\begin{example}
Suppose $(\A,\BE,\fs)$ is an
abelian category. 
Then the notions defined in \cref{def:E-s-notions} are equivalent to the classical ones for $\A$. That is, an object in $\A$ is simple if and only if it is $(\BE,\fs)$-simple, etc.
\end{example}

\begin{example}
Suppose $(\A,\BE,\fs)$ is an
exact category $(\A,\mathscr{E})$. 
Then the notions in \cref{def:E-s-notions}\ref{item:E-s-subobject}--\ref{item:E-s-simple} agree with those defined in 
\cite[Defs.\ 3.1, 3.3]{BHLR}. 
\end{example}

Since we cannot formally call a zero object in an additive category simple, 
we make the following definition.

\begin{definition}
\label{def:simplie-like-atom-like}
We say that a zero object $0\in\A$ is \emph{$(\BE, \fs)$-simple-like} if, whenever there is an extriangle 
    $
    \begin{tikzcd}[column sep=0.5cm,cramped]
        A\arrow{r}{a}
        &0\arrow{r}
        &C\arrow[dashed]{r}{}
        &{},
    \end{tikzcd}
    $ 
    the $(\BE,\fs)$-inflation 
    $a$ is an isomorphism (equivalently, $A= 0$).
\end{definition}

Asking that $0$ is $(\BE,\fs)$-simple-like is necessary for a length extriangulated category to have the Jordan-H{\"{o}}lder property (see \cref{def:E-s-notions}). 
Indeed, 
suppose that there is an extriangle 
$
\label{equation:X-triangle}
\begin{tikzcd}[column sep=0.5cm,ampersand replacement=\&,cramped]
X \arrow{r}\& 0 \arrow{r}\& Y \arrow[dashed]{r}\& {}
\end{tikzcd}
$
in $(\A,\BE,\fs)$ 
with $X \not\cong 0$. 
For example, take $(\A,\BE,\fs)$ to be a non-zero triangulated category (or an extension-closed subcategory thereof large enough to contain $X$ and $X[1]$), so that $Y \iso X[1]$. One obtains 
arbitrarily long 
filtrations 
\[
\begin{tikzcd}[row sep=0.5cm]
\cdots \arrow{rr}&&0 \arrow{rr}&& X \arrow{rr} \arrow{dl}&& 0 \arrow{rr}\arrow{dl}&& X \arrow{dl}\\
&&&X \arrow[dashed]{ul}&& Y \arrow[dashed]{ul}&& X \arrow[dashed]{ul}
\end{tikzcd}
\] of $X$ by non-isomorphisms, 
which is 
an obstruction to the Jordan-H{\"{o}}lder property.

The following observation is due to Carlo Klapproth; it identifies a sufficiency condition for the zero object to be simple-like.

\begin{prop}
\label{prop:A-has-simple-then-zero-simple}
    If $(\A,\BE,\fs)$ has an $(\BE,\fs)$-simple object $S$, then $0$ must be $(\BE,\fs)$-simple-like.
\end{prop}
\begin{proof}
    Assume there is an extriangle 
    $
    \begin{tikzcd}[column sep=0.5cm,cramped]
        A\arrow{r}{a}
        &0\arrow{r}
        &C\arrow[dashed]{r}{\delta}
        &{}
    \end{tikzcd}
    $
    in $\A$. Since 
    $
    \begin{tikzcd}[column sep=0.7cm,cramped]
        S\arrow{r}{\id{S}}
        &S\arrow{r}
        &0 \arrow[dashed]{r}{{}_{S}0_{0}}
        &{}
    \end{tikzcd}$
    is also an extriangle, we obtain the direct sum
    \[
    \begin{tikzcd}[ampersand replacement=\&, column sep=1.5cm]
        A \oplus S \arrow{r}{\begingroup
        \renewcommand\thickspace{\kern5pt}
        \begin{psmallmatrix}
        0 & \id{S}
        \end{psmallmatrix}
        \endgroup}
        \& S\arrow{r}
        \& C \oplus 0\arrow[dashed]{r}{\delta \oplus {}_{S}0_{0}}
        \&{.}
    \end{tikzcd}
    \]
    As $S$ is $(\BE,\fs)$-simple, the morphism 
    $
    \begingroup
        \renewcommand\thickspace{\kern5pt}
        \begin{psmallmatrix}
        0 & \id{S}
        \end{psmallmatrix}
    \endgroup
    $
    is either zero or an isomorphism. 
    It cannot be zero because $S$ is non-zero, and hence we have an isomorphism 
    $
    \begingroup
        \renewcommand\thickspace{\kern5pt}
        \begin{psmallmatrix}
        0 & \id{S}
        \end{psmallmatrix}
    \endgroup
    \colon 
    A\oplus S \to S
    $.
    This implies there is a morphism 
    $
    \begin{psmallmatrix}
        x \\ y
    \end{psmallmatrix}
    \colon S \to A\oplus S
    $, such that 
    $
    \begin{psmallmatrix}
        \id{A} & 0 \\
        0 & \id{S}
    \end{psmallmatrix}
        = \begin{psmallmatrix}
        x \\ y
    \end{psmallmatrix} \begingroup
        \renewcommand\thickspace{\kern5pt}
        \begin{psmallmatrix}
        0 & \id{S}
        \end{psmallmatrix}
    \endgroup
        = \begin{psmallmatrix}
        0 & x \\
        0 & y
        \end{psmallmatrix}
    $.
    Thus, $\id{A} = 0$ and hence 
    $A= 0$. 
    Therefore, $0$ is $(\BE,\fs)$-simple-like.
\end{proof}

The zero object being $(\BE,\fs)$-simple-like has some nice consequences on the $(\BE,\fs)$-simple objects in $\A$ and on the Grothendieck monoid. 
First, we see that $(\BE,\fs)$-simples satisfy a stronger property that  resembles the classical definition in an abelian category.

\begin{lem}
\label{lem:0-simple-like-implies-stronger-simple}
Suppose $(\A,\BE,\fs)$ is weakly idempotent complete. 
The following are equivalent for a non-zero object $S\in\A$.
\begin{enumerate}[label=\textup{(\roman*)}]

    \item\label{item:B-simple} $S$ is $(\BE,\fs)$-simple.

    \item\label{item:B-strongly-simple} For each extriangle 
    $\begin{tikzcd}[column sep=0.5cm,cramped]
    A \arrow{r}{a}& S \arrow{r}& C \arrow[dashed]{r}& {,}
    \end{tikzcd}$
    we have either $A = 0$ or $a$ is an isomorphism.
    
    \item\label{item:B-WWZZ-simple} 
    There does not exist an extriangle 
    $\begin{tikzcd}[column sep=0.5cm,cramped]
    A \arrow{r}{a}& S \arrow{r}& C \arrow[dashed]{r}& {}
    \end{tikzcd}$
    with both $A$ and $C$ non-zero. That is, $S$ is `simple' in the sense of \cite{WWZZ}.
\end{enumerate}
\end{lem}
\begin{proof}
It is clear that \ref{item:B-strongly-simple} implies 
\ref{item:B-simple}.

\ref{item:B-simple} $\Rightarrow$ \ref{item:B-strongly-simple}:\; 
Suppose $S\neq 0$ is $(\BE,\fs)$-simple and let 
\begin{equation}
\label{eqn:B-simple}
\begin{tikzcd}[column sep=0.5cm,cramped]
    A \arrow{r}{a}& S \arrow{r}& C \arrow[dashed]{r}& {}
\end{tikzcd}
\end{equation}
be an extriangle in $(\A,\BE,\fs)$. 
If $a$ is an isomorphism, then we are done, so suppose that $a=0$. 
From the dual of \cref{lem:decomposable}, we deduce that \eqref{eqn:B-simple} is isomorphic to the direct sum of 
two extriangles, with one of the form 
$\begin{tikzcd}[column sep=0.5cm,cramped]
    0 \arrow{r}{}& S \arrow{r}{\id{S}}& S \arrow[dashed]{r}& {}
\end{tikzcd}$
and the other of the form 
$\begin{tikzcd}[column sep=0.5cm,cramped]
    A \arrow{r}{}& 0 \arrow{r}& C' \arrow[dashed]{r}& {.}
\end{tikzcd}$
Since $0$ is $(\BE,\fs)$-simple-like by \cref{prop:A-has-simple-then-zero-simple}, we must have that $A=0$. 
%Thus, we have \ref{item:B-simple}$\Leftrightarrow$\ref{item:B-strongly-simple}.

\ref{item:B-strongly-simple} $\Rightarrow$ \ref{item:B-WWZZ-simple} :\; 
If there were an extriangle 
$\begin{tikzcd}[column sep=0.5cm,cramped]
    A \arrow{r}{a}& S \arrow{r}& C \arrow[dashed]{r}& {}
\end{tikzcd}$
with both $A$ and $C$ non-zero, 
then $a$ would be an isomorphism by \ref{item:B-strongly-simple} and $C$ would be zero, a contradiction.

\ref{item:B-WWZZ-simple} $\Rightarrow$ \ref{item:B-strongly-simple}:\; 
Suppose 
$\begin{tikzcd}[column sep=0.5cm,cramped]
    A \arrow{r}{a}& S \arrow{r}& C \arrow[dashed]{r}& {}
\end{tikzcd}$
is an extriangle. By \ref{item:B-WWZZ-simple}, either $A=0$ or $C=0$. In the former, we are done immediately. The latter implies $a$ is an isomorphism, and we are also done.
\end{proof}

\begin{cor}
\label{cor:simple-isomorphic-to-X}
Suppose $(\A,\BE,\fs)$ is weakly idempotent complete. 
Let $\X$ be a class of non-zero objects of $\A$. 
Then each 
$(\restr{\BE}{\F(\X)},\restr{\fs}{\F(\X)})$-simple object $M\in\F(\X)$ is isomorphic to an object in $\X$. 
\end{cor}
\begin{proof}
If \eqref{eqn:X-filtration} is an $\X$-filtration of $M\in\A$, then we have $t\geq 1$ and that 
there is an extriangle
$
\begin{tikzcd}[column sep=0.5cm,cramped]
    M_{t-1} \arrow{r}{a}& M \arrow{r}{}& X_{t} \arrow[dashed]{r}{\xi_{t}} &{}
\end{tikzcd}
$
in  $\F(\X)$
with $X_{t}\in\X$. 
We deduce that $M_{t-1} = 0$ 
or that 
$a$ 
is an isomorphism by \cref{lem:0-simple-like-implies-stronger-simple}. 
The latter implies that $X_{t}$ is a zero object, which is contrary to our assumption on $\X$. 
Thus, we see that $M_{t-1} = 0$ and hence $M\iso X_{t}$. 
\end{proof}

The converse to \cref{cor:simple-isomorphic-to-X} fails in general as $\X$ can contain redundancy, in the sense that it may contain the direct sum of two other objects already in $\X$.

\begin{prop}
\label{prop:0-simple-like-implies-simples-have-comp-length-1}
Suppose $(\A,\BE,\fs)$ is weakly idempotent complete. 
The following are equivalent for an object $M\in\A$.
\begin{enumerate}[label=\textup{(\roman*)}]
    \item\label{item:M-simple} $M$ is $(\BE,\fs)$-simple.
    \item\label{item:M-has-comp-series-only-of-length-1} $M$ has an $(\BE,\fs)$-composition series and every $(\BE,\fs)$-composition series of $M$ has \mbox{length $1$.}
\end{enumerate}
\end{prop}

\begin{proof}
\ref{item:M-simple} $\Rightarrow$ \ref{item:M-has-comp-series-only-of-length-1}:\; 
Suppose $\xi = (\xi_{0},\xi_{1},\ldots,\xi_{t})$ from \eqref{eqn:X-filtration} is an $(\BE,\fs)$-composition series of $M$. 
If $t=0$, then $M=0$ which is a contradiction. If $t=1$, then we are done. 
Now assume for contradiction that $t\geq 2$. 
From the $(\BE,\fs)$-composition series 
% $\xi = (\xi_{0},\xi_{1},\ldots,\xi_{t})$, 
$\xi$, 
we see that there is an extriangle
$\begin{tikzcd}[column sep=0.5cm,cramped]
M_{t-1} \arrow{r}{a}
    & M \arrow{r}{b}
    & X_{t} \arrow[dashed]{r}{\xi_{t}}
    &{}
\end{tikzcd}$
with $X_{t}$ an  $(\BE,\fs)$-simple.
Since we are supposing $M$ is $(\BE,\fs)$-simple, 
we have that $M_{t-1}$ is zero or $a$ is an isomorphism by \cref{lem:0-simple-like-implies-stronger-simple}. 
The latter would imply $X_{t} = 0$, which is impossible, so we must have $M_{t-1}=0$. 
In particular, we see that we have an extriangle 
$\begin{tikzcd}[column sep=0.8cm,cramped]
M_{t-2} \arrow{r}{}
    & 0 \arrow{r}{b}
    & X_{t-1} \arrow[dashed]{r}{\xi_{t-1}}
    &{,}
\end{tikzcd}$
so $M_{t-2}$ is forced to be the zero object as $0$ is $(\BE,\fs)$-simple-like by \cref{prop:A-has-simple-then-zero-simple}. This implies that $X_{t-1} = 0$, which is a contradiction so $t\not\geq 2$.

\ref{item:M-has-comp-series-only-of-length-1} $\Rightarrow$ \ref{item:M-simple}:\; 
By assumption, $M$ has an $(\BE,\fs)$-composition series of length $1$. That is, there is an $\Sim(\A,\BE,\fs)$-filtration of the form
\[
\begin{tikzcd}[column sep=0.7cm, row sep=0.5cm]
0 
% = M_{-1} 
\arrow{rr}&&
0 
% = M_{0} 
\arrow{rr}\arrow{dl}&& M_{1} = M \arrow{dl}\\
& 0
% = X_{1}
\arrow[dashed]{ul}{\xi_{0}}&&X_{1}\arrow[dashed]{ul}{\xi_{1}}&
\end{tikzcd}
\]
in which $X_{1}$ is $(\BE,\fs)$-simple. In particular, $M \iso X_{1}$ and $M$ is $(\BE,\fs)$-simple.
\end{proof}

%%%%%%%%%%%%%%%%%%%%%%%%%%%%%%%%%%%%%%
%%%%%%%%%%%%%%%%%%%%%%%%%%%%%%%%%%%%%%

\subsection{The Grothendieck monoid of an extriangulated category}
\label{subsec:grothendieck-monoid}
For the remainder of this section, we work under the following setup.

\begin{setup}
\label{setup:skel-small-extriangulated-cat}
By $(\A,\BE,\fs)$ we denote a skeletally small, extriangulated category. 
\end{setup}

The Grothendieck group of a certain extriangulated category was first considered in \cite{PPPP}. 
The \emph{Grothendieck group} $K_{0}(\A,\BE,\fs)$ of $(\A,\BE,\fs)$ is defined analogously to the case of exact categories or triangulated categories \cite[Sec.~4]{ZZ21} (also \cite[Def.\ 4.1]{Haugland-the-grothendieck-group-of-an-n-exangulated-category}): 
let $F(\A)$ denote the free abelian group generated by the isomorphism classes $[X]$ of all objects $X \in \A$, then the Grothendieck group $K_{0}(\A,\BE,\fs)$ is defined as the quotient
of $F(\A)$ by the subgroup generated by 
each relation 
$
[X] + [Z] - [Y] 
$ 
whenever 
$
\begin{tikzcd}[column sep=0.5cm,cramped]
X\arrow{r}&Y\arrow{r}&Z\arrow[dashed]{r}&{}
\end{tikzcd}
$
is an extriangle.

Although the Grothendieck group is a suitable categorical invariant for 
e.g.\ skeletally small, abelian, length categories, 
it is not quite sensitive enough for arbitrary skeletally small, exact categories. 
Indeed, the methods employed in \cite{BeGr} and \cite{E19} to study filtrations in an exact category use the Grothendieck monoid associated to the exact structure. 
We explain the corresponding idea for extriangulated categories. 
These concepts were formulated independently in \cite[Sec.\ 2.2]{EnomotoSaito} in order to, amongst other things, classify certain (e.g.\ Serre) subcategories of an extriangulated category.

All monoids considered here are commutative, with identity $0$ and additive notation. 
We denote 
the set of isomorphism classes of objects of $\A$ by $\Iso \A$. 
We say a function $f\colon  \Iso \A \to \Gamma$ to some monoid $\Gamma$ is \emph{additive on extriangles} if $f[0]=0$ and $f[Y]=f[X]+f[Z]$ for every extriangle
$
\begin{tikzcd}[column sep=0.5cm,cramped]
X\arrow{r}&Y\arrow{r}&Z\arrow[dashed]{r}&{}
\end{tikzcd}
$ 
(see also \cite[Def.\ 3.3]{JorgensenShah-the-index-with-respect-to-a-rigid-subcategory}).

\begin{definition}
\label{def:grothendieck-monoid}
The \emph{Grothendieck monoid} 
of a skeletally small, extriangulated category $(\A,\BE,\fs)$ is a monoid $\M(\A,\BE,\fs)$ together with a map $q \colon \Iso \A \to \M(\A,\BE,\fs)$ which is additive on extriangles, such that: 
every function $f\colon \Iso\A \to \Gamma$ to a monoid $\Gamma$ that is additive on extriangles factors uniquely through $q$. 
\end{definition}

The existence of $\M(\A,\BE,\fs)$ can be shown like in \cite[Prop.\ 3.3]{E19} or \cite[Sec.\ 2.3]{BeGr} as follows. 
First, the set $\Iso \A$ becomes a monoid under the operation 
$
[X] + [Y] = [X\oplus Y]
$. 
Note that this implies that $[0]$ is the identity element of the monoid. 
Then, by \cite[p.\ 14]{Grillet-semigroups}, there is a \emph{monoid congruence} (see \cite[Def.\ A.5]{E19}), denoted  $\approx$, on $\Iso\A$ generated by the binary relation $\sim$, where we have 
$
[Y] \sim [X] + [Z]
$
whenever 
$
\begin{tikzcd}[column sep=0.5cm,cramped]
X\arrow{r}&Y\arrow{r}&Z\arrow[dashed]{r}&{}
\end{tikzcd}
$
is an extriangle in $(\A,\BE,\fs)$. 
It is then clear that 
$
\M(\A,\BE,\fs) 
    \deff \Iso\A / \approx
$ 
is a monoid with the properties of 
\cref{def:grothendieck-monoid}. 
Note that every element $a\in\M(\A,\BE,\fs)$ is a finite sum of elements from $\Iso\A$, and so the monoid operation implies $a = [A]$ for some $A\in\A$. 
Furthermore, the Grothendieck group $K_{0}(\A,\BE,\fs)$ of 
$(\A,\BE,\fs)$ can be realized as the \emph{group completion} of the Grothendieck monoid $\M(\A,\BE,\fs)$ (see \cite[Def.\ A.10, Prop.\ 3.9]{E19}). 
We refer the reader to \cite{BeGr,E19} for several properties of the Grothendieck monoid for exact categories.

\begin{definition} 
\label{def:atom-reduced}
Let $\M$ be a monoid. 
\begin{enumerate}[label=\textup{(\roman*)}]
    \item $\M$ is \emph{reduced} if $x +y =0$ implies $x =y = 0$ in $\M$ for all $x, y \in \M$.
    \item A non-zero element $0\neq a \in \M$ is an \emph{atom} if $a = x+y$ implies that $x=0 $ or $y=0$ for all  $x, y \in \M$. By $\Atom(\M)$ we denote the set of all atoms of $\M$.
\end{enumerate}
\end{definition}

In order to relate 
\cref{def:simplie-like-atom-like}
to whether the monoid $\M(\A,\BE,\fs)$ is reduced, we will use the $\X$-filtration defined in \cref{subsec:filtered-objects-in-an-extriangulated-category} with $\X = \A$. 
These are the analogues of the `inflation series' introduced in \cite[p.\ 11]{E19}. 
If 
$\begin{tikzcd}[column sep=0.5cm,cramped]
A \arrow{r}{a}& B \arrow{r}& C \arrow[dashed]{r}&{}
\end{tikzcd}$
is an extriangle, then we will use the notation 
$\begin{tikzcd}[column sep=0.5cm,cramped]
A \arrow[tail]{r}{a}& B
\end{tikzcd}$ 
to indicate that $a$ is an inflation.

\begin{definition}
\label{def:inflation-series-defs}
Let $A,B,M,M'$ be objects of $\A$. 
\begin{enumerate}[label=\textup{(\roman*)}]

    \item By an \emph{$(\BE,\fs)$-inflation series from $A$ to $B$}, we mean a sequence of inflations of the form
\begin{equation}\label{eqn:inflation-series-from-A-to-B}
\begin{tikzcd}%[column sep=0.5cm]
A=M_{0} \arrow[tail]{r}{f_{0}} & M_{1}\arrow[tail]{r}{f_{1}} & \cdots \arrow[tail]{r}{f_{t-2}} & M_{t-1} \arrow[tail]{r}{f_{t-1}} & M_{t} = B,
\end{tikzcd}
\end{equation}
where, for each $i = 1, \ldots, t$, there is an extriangle
$
\begin{tikzcd}[cramped]%[column sep=0.5cm]
M_{i-1} \arrow{r}{f_{i-1}} & M_{i} \arrow{r} & X_{i}
\arrow[dashed]{r}{\xi_{i}} & {}
\end{tikzcd}
$
with $X_{i} \neq 0$. 
We suggestively write 
$M_{i}/M_{i-1}$
for the object 
$X_{i}$.

    \item If $A=0$, then an $(\BE,\fs)$-inflation series from $A$ to $B$ is called more simply an \emph{$(\BE,\fs)$-inflation series of $B$}. In particular, $(\BE,\fs)$-inflation series of $B$ correspond to $\A$-filtrations 
    of $B$ in the sense of \cref{def:X-filtration} 
    with the extra condition that the factors 
    $X_{i}$ are non-zero for all $1\leq i \leq t$.  
    
    \item Suppose 
     $
    \begin{tikzcd}[column sep=0.5cm,cramped]
    0=M_{0} \arrow[tail]{r}{} 
    % & M_{1}\arrow[tail]{r}{}
    & \cdots \arrow[tail]{r}{} 
    % & M_{t-1} \arrow[tail]{r}{} 
    & M_{t} = M
    \end{tikzcd}
    $
    and 
     $
    \begin{tikzcd}[column sep=0.5cm,cramped]
    0=M'_{0} \arrow[tail]{r}{} 
    % & M'_{1}\arrow[tail]{r}{} 
    & \cdots \arrow[tail]{r}{} 
    % & M'_{s-1} \arrow[tail]{r}{} 
    & M'_{s} = M'
    \end{tikzcd}
    $
    are $(\BE,\fs)$-inflation series of $M$ and $M'$, respectively.
    We call them \emph{isomorphic} if $s=t$ and there is a permutation $\sigma\in\Sym{n}$ such that 
    $M_{i}/M_{i-1}$ and $M'_{\sigma(i)}/M'_{\sigma(i-1)}$ are isomorphic for all $1\leq i\leq t$.
\end{enumerate}
\end{definition}

With this terminology in place, one can check the proof of \cite[Prop.\ 3.4]{E19} carries over to our setting, which yields the following. 

\begin{prop}
\label{prop:e19-prop3-4-analogue}
    For objects $M,M'\in\A$, we have $[M]=[M']$ in $\M(\A,\BE,\fs)$ if and only if there exists a sequence of objects 
    $
    M = A_{0}, A_{1}, \ldots, A_{n} = M'
    $ 
    in $\A$ such that $A_{i-1}$ and $A_{i}$ have isomorphic $(\BE,\fs)$-inflation series for each $1\leq i\leq n$. 
\end{prop}

\begin{remark}
We note that `having isomorphic $(\BE,\fs)$-inflation series' is not a transitive relation. 
That is, even if
$
M = A_{0}, A_{1}, \ldots, A_{n} = M'
$ 
are objects of $\A$ with $A_{i-1}$ and $A_{i}$ having isomorphic $(\BE,\fs)$-inflation series for each $1\leq i\leq n$, this does not a priori 
imply $M$ and $M'$ have isomorphic $(\BE,\fs)$-inflation series. 
Indeed, suppose $(\A,\BE,\fs)$ is a triangulated with objects $X,Y$ such that $X\niso Y$. Then there are triangles 
$
\begin{tikzcd}[column sep=0.5cm,cramped]
X \arrow{r}& 0 \arrow{r}& X[1] \arrow{r}{\id{X[1]}}&[0.5cm] X[1]
\end{tikzcd}
$, 
$
\begin{tikzcd}[column sep=0.5cm,cramped]
0 \arrow{r}& X \arrow{r}{\id{X}}&[0.3cm] X \arrow{r}& 0
\end{tikzcd}
$ 
and 
$
\begin{tikzcd}[column sep=0.5cm,cramped]
X \arrow{r}& X\oplus X[1] \arrow{r}& X[1] \arrow{r}{0}& X[1]
\end{tikzcd}
$, 
and similar ones involving $Y$. 
These triangles can be used to show that $0$ and $X\oplus X[1]$ have isomorphic $(\BE,\fs)$-inflation series, and similarly that $0$ and $Y\oplus Y[1]$ have isomorphic $(\BE,\fs)$-inflation series. 
However, one cannot a priori deduce that $X\oplus X[1]$ and $Y\oplus Y[1]$ have isomorphic $(\BE,\fs)$-inflation series.
\end{remark}

\begin{prop}
\label{prop:consequences-of-zero-simple-like}
Suppose $0$ is $(\BE,\fs)$-simple-like. 
Then the following statements hold. 
\begin{enumerate}[label=\textup{(\roman*)}]

    \item\label{item:0-simple-like-then-0-in-monoid-implies-iso-to-0} $[M]=[0]$ in $\M(\A,\BE,\fs)$ if and only if $M\iso 0$ in $\A$.
    
    \item\label{item:M-reduced} $\M(\A,\BE,\fs)$ is reduced.
    
    \item\label{item:bijection-simples-atoms} 
    If $(\A,\BE,\fs)$ is weakly idempotent complete, then there is a bijection as follows.
    \begin{align*}
    \Set{[S]\in\Iso\A | S\in\Sim(\A,\BE,\fs)}
& \longrightarrow 
    \Atom( \M(\A,\BE,\fs) )\\[5pt]
    [S]
    &
    \xmapsto{\hspace{11pt}}
    [S]
\end{align*}
\end{enumerate}
\end{prop}

\begin{proof}
\ref{item:0-simple-like-then-0-in-monoid-implies-iso-to-0}\;
If $M\iso 0$, then $[M]=[0]$ in $\Iso\A$ and hence also in $\M(\A,\BE,\fs)$. 
Conversely, assume
$[M]=[0]$ in $\M(\A,\BE,\fs)$. 
By \cref{prop:e19-prop3-4-analogue} there is a sequence of objects 
$
0 = A_{0}, A_{1}, \ldots, A_{n} = M
$ 
in $\A$ such that $A_{i-1}$ and $A_{i}$ have isomorphic $(\BE,\fs)$-inflation series for each $1\leq i\leq n$. 
Since $0$ is $(\BE,\fs)$-simple-like, \emph{any} $(\BE,\fs)$-inflation series 
$
\begin{tikzcd}[column sep=0.5cm,cramped]
    0=M_{0} \arrow[tail]{r}{} 
    % & M_{1}\arrow[tail]{r}{}
    & \cdots \arrow[tail]{r}{} 
    % & M_{t-1} \arrow[tail]{r}{} 
    & M_{t} = 0
\end{tikzcd}
$
of $0$ has $M_{i}/M_{i-1} = 0$ for all $1\leq i\leq t$. 
Since $A_{1}$ has an $(\BE,\fs)$-inflation series isomorphic to one of $0$, it has an $(\BE,\fs)$-inflation series 
$
\begin{tikzcd}[column sep=0.5cm,cramped]
    0=M^{1}_{0} \arrow[tail]{r}{} 
    % & M_{1}\arrow[tail]{r}{}
    & \cdots \arrow[tail]{r}{} 
    % & M_{t-1} \arrow[tail]{r}{} 
    & M^{1}_{t_{1}} = A_{1}
\end{tikzcd}
$
in which $M^{1}_{i}/M^{1}_{i-1} = 0$ for all $1\leq i\leq t_{1}$. This implies, 
$0 = M^{1}_{0}\iso M^{1}_{1}\iso \cdots \iso M^{1}_{t_{1}} = A_{1}$. 
Continuing in this way, we see that $M = A_{n} \iso 0$.

\ref{item:M-reduced}\;
Assume $[X] + [Y] = [0]$ for some $X,Y\in\A$. 
Thus, $[X\oplus Y] = [0]$, 
so by part \ref{item:0-simple-like-then-0-in-monoid-implies-iso-to-0} we have $X\oplus Y\iso 0$.
Thus, $X\iso 0\iso Y$, whence $[X]=[0]$ and $[Y]=[0]$ by \ref{item:0-simple-like-then-0-in-monoid-implies-iso-to-0}.

\ref{item:bijection-simples-atoms}\;
In light of 
\cref{prop:e19-prop3-4-analogue}, 
\cref{lem:0-simple-like-implies-stronger-simple} 
and 
part \ref{item:0-simple-like-then-0-in-monoid-implies-iso-to-0}, 
one may carefully recycle the argument of \cite[Prop.\ 3.6]{E19} to our setting. See also \cite[Lem.\ 5.3]{BeGr}.
\end{proof}

%%%%%%%%%%%%%%%%%%%%%%%%%%%%%%%%%%%%%%%%%%%%%%%%%%
%%%%%%%%%%%%%%%%%%%%%%%%%%%%%%%%%%%%%%%%%%%%%%%%%%

\subsection{A characterisation of length Jordan-H{\"{o}}lder extriangulated categories}
\label{subsec:characterisation-of-JHP}

We still assume \Cref{setup:skel-small-extriangulated-cat}. 
We need one last lemma and definition before proving our main result (see \cref{thm:JH-iff-free-monoid-iff-free-group}) of this section.
The lemma follows from the proof of \cite[Lem.~2.7]{WWZ21}.

\begin{lem}
\label{lem:comp-series-of-A-and-C-gives-one-for-B}
Suppose  
$\begin{tikzcd}[column sep=0.5cm,cramped]
A \arrow{r}{f} & B \arrow{r}{g} & C \arrow[dashed]{r}{\eta} & {}
\end{tikzcd}$
is an extriangle 
in $(\A,\BE,\fs)$. 
For each $(\BE,\fs)$-inflation series of $C$, 
there exists an $(\BE,\fs)$-inflation series \eqref{eqn:inflation-series-from-A-to-B} from $A$ to $B$ with the same factors up to isomorphism such that $f_{t-1}f_{t-2}\cdots f_{0} = f$. 
In particular, given $\X$-filtrations of $A$ and $C$ of lengths $u$ and $v$, respectively, one may construct an $\X$-filtration of $B$ of length $u+v$ with the same $\X$-factors up to isomorphism. 
\end{lem}

\begin{definition}
\cite[Def.\ A.16]{E19}
Let $\M$ be a monoid. A monoid homomorphism $\nu\colon \M\to \bN$ is called \emph{length-like} if $\nu(x) = 0$ implies that $x=0$ for all $x\in\M$. 
\end{definition}

The next result and its proof are heavily inspired by \cite[Thm.\ 4.9]{E19}.

\begin{thm}
\label{thm:JH-iff-free-monoid-iff-free-group}
Let $(\A,\BE,\fs)$ be a skeletally small, weakly idempotent complete extriangulated category. 
Then the following conditions are equivalent.
\begin{enumerate}[label=\textup{(\roman*)}]
    \item The extriangulated category $(\A,\BE,\fs)$ is length and Jordan-H{\"{o}}lder. 
    \item 
    The Grothendieck monoid $\M(\A,\BE,\fs)$ is free with basis the image of $\Sim(\A,\BE,\fs)$.
    \item 
    The Grothendieck group $K_{0}(\A,\BE,\fs)$ is free with basis the image of $\Sim(\A,\BE,\fs)$.
\end{enumerate}
\end{thm}

\begin{proof}
In light of \cref{prop:consequences-of-zero-simple-like}, the equivalence of (ii) and (iii) follows from \cite[Thm.\ A.20]{E19}. Thus, it suffices to prove (i) $\Leftrightarrow$ (ii). 
Furthermore, if $\A$ has no $(\BE,\fs)$-simple objects, then we observe that (i) and (ii) are trivially equivalent because it forces $\A$ to be zero in these equivalent case. Thus, we also assume $\A$ has at least one $(\BE,\fs)$-simple object in the following. In particular, by \cref{prop:A-has-simple-then-zero-simple}, we know that $0$ is $(\BE,\fs)$-simple-like.

(i) $\Rightarrow$ (ii):\; 
By \cref{prop:consequences-of-zero-simple-like}\ref{item:bijection-simples-atoms}, we have 
$\Atom( \M(\A,\BE,\fs) ) = \set{[S] | S\in\Sim(\A,\BE,\fs) }$. 
The remainder of the argument is as in the proof of $(1)\Rightarrow (2)$ of \cite[Thm.\ 4.9]{E19}, using \cref{lem:comp-series-of-A-and-C-gives-one-for-B} in place of \cite[Prop.\ 2.5(2)]{E19}.

(ii) $\Rightarrow$ (i):\; 
By \cref{prop:consequences-of-zero-simple-like}, we have that $\M(\A,\BE,\fs)$ is reduced 
and again the equality $\Atom( \M(\A,\BE,\fs) ) = \set{[S] | S\in\Sim(\A,\BE,\fs) }$. 
Since $\M(\A,\BE,\fs)$ is free, and hence factorial \cite[Prop.\ A.9]{E19},  we have that there is a length-like function $\nu\colon \M(\A,\BE,\fs) \to \bN$ satisfying $\nu([S]) = 1$ for each $[S]\in\Atom( \M(\A,\BE,\fs) )$ by \cite[Lem.\ A.18]{E19}.

To show each $M\in\A$ admits an $(\BE,\fs)$-composition series, we induct on $\nu([M])\in\bN$. 
Since $\nu$ is length-like and $0$ is $(\BE,\fs)$-simple-like, we have that $\nu([M]) = 0$ if and only if $[M] = 0$ in $\M(\A,\BE,\fs)$. In this case, $M$ has an $(\BE,\fs)$-composition series of length $0$. 
Assume $\nu([M]) = 1$. We claim that this implies $M$ is $(\BE,\fs)$-simple. 
To this end, suppose there is an extriangle 
$
    \begin{tikzcd}[column sep=0.5cm,cramped]
        A\arrow{r}{a}
        &M\arrow{r}
        &C\arrow[dashed]{r}{}
        &{}.
    \end{tikzcd}
$ 
Then $1 = \nu([M]) = \nu([A]) + \nu([C])$ since $\nu$ is a monoid homomorphism. 
Therefore, $\nu([A]) = 0$ or $\nu([C]) = 0$, which implies $A\iso 0$ or $C\iso 0$, and hence that $M$ is $(\BE,\fs)$-simple. 
Now suppose $\nu([M]) >1$ and that the claim holds for objects $Y$ with $\nu([Y])<\nu([M])$. 
If $M$ were $(\BE,\fs)$-simple, then we would have that $[M]$ is an atom and $\nu([M]) = 1$, which is not possible. So there exists an extriangle 
$
    \begin{tikzcd}[column sep=0.5cm,cramped]
        A\arrow{r}
        &M\arrow{r}
        &C\arrow[dashed]{r}{}
        &{}
    \end{tikzcd}
$ 
with $A$ and $C$ both non-zero and, in particular, $\nu([A]),\nu([C])\geq 1$.
Since $\nu([M]) = \nu([A]) + \nu([C])$ and $\nu$ takes values in $\bN$, 
we see that $\nu([A]),\nu([C]) < \nu([M])$. Hence, $A$ and $C$ both admit an $(\BE,\fs)$-composition series. Then we may apply \cref{lem:comp-series-of-A-and-C-gives-one-for-B} (with $\X = \Sim(\A,\BE,\fs)$) to obtain an $(\BE,\fs)$-composition series of $M$.

Now we show $(\A,\BE,\fs)$ is a length extriangulated category. 
So suppose \eqref{eqn:X-filtration} is an $(\BE,\fs)$-composition series of $M$. 
Then $\nu([M]) = \nu([X_{1}]) + \cdots + \nu([X_{t}]) = t$ as each $X_{i}$ is $(\BE,\fs)$-simple. 
In particular, each $(\BE,\fs)$-composition series of $M$ is of finite length $\nu([M])\in\bN$. 

To see that $(\A,\BE,\fs)$ is Jordan-H{\"{o}}lder, let 
$(\xi_{1},\ldots,\xi_{t})$ 
and 
$(\xi'_{1},\ldots,\xi'_{s})$ 
be two $(\BE,\fs)$-composition series 
of $M\in\A$. 
Since $\M(\A,\BE,\fs)$ is a free monoid, we must have that 
$s=t$ and there is a permutation $\sigma\in\Sym{n}$ such that 
$[X_{i}] = [X'_{\sigma(i)}]$. 
It follows that $X_{i} \iso X'_{\sigma(i)}$ by \cref{prop:consequences-of-zero-simple-like}\ref{item:bijection-simples-atoms}. 
\end{proof}

%%%%%%%%%%%%%%%%%%%%%%%%%%%%%%%%%%%%%%%%%%%%%%%%%%%%%%%%
%%%%%%%%%%%%%%%%%%%%%%%%%%%%%%%%%%%%%%%%%%%%%%%%%%%%%%%%

\section{Stratifying systems}
\label{sec:stratifying-systems}

In this section, we define stratifying systems in the extriangulated setting (see \cref{def:stratifying-system}). We show that under mild assumptions, the extriangulated category $\F(\Phi)$ of filtered objects induced by an $\BE$-stratifying system $\Phi$ has the Jordan-H{\"{o}}lder property (see \cref{thm:JHP-for-Theta-filtrations}). 
Along the way, we also show that each $\BE$-stratifying system 
can be completed to a projective $\BE$-stratifying system, and hence $\F(\Phi)$ has enough relative projectives (see \cref{thm:existence-of-proj-ESS}), laying the groundwork for future study. 
For examples of the notions and results we present in this section, we refer the reader to \cref{sec:examples}.

Suppose $\ring$ is a commutative artinian ring. 
We say that the extriangulated category $(\A,\BE,\fs)$ is \emph{$\ring$-linear} if $\A$ is $\ring$-linear in the usual sense and, moreover, if $\BE(C,A)$ is an $\ring$-module for all $A,C\in\A$ such that $\BE$ is $\ring$-bilinear. 
If $(\A,\BE,\fs)$ is $\ring$-linear, then we say it is \emph{$\BE$-finite} if $\BE(C,A)$ is a finite length $\ring$-module for all $A,C\in\A$.

\begin{definition}
\label{def:artin-extriangulated-category}
\cite[Def.\ 2.10]{Zh19} 
Suppose $\ring$ is a commutative artinian ring. 
An $R$-linear extriangulated category $(\A,\BE,\fs)$ is called an 
\emph{artin $\ring$-linear extriangulated category} 
if $\A$ is idempotent complete and $\Hom$-finite (over $\ring$), and $(\A,\BE,\fs)$ is $\BE$-finite. 
\end{definition}

By \cite[Thm.~6.1]{Shah-KRS-decompositions-in-Hom-finite-additive-categories}, an artin $\ring$-linear extriangulated category is Krull-Schmidt (and hence also weakly idempotent complete, see \cref{sec:extriangulated-categories}). 
For the remainder of \cref{sec:stratifying-systems}, we work in the following context.

\begin{setup}
\label{setup:A-is-artin-extriangulated}
We assume $\ring$ is a commutative artinian ring and that $(\A,\BE,\fs)$ is an artin $\ring$-linear extriangulated  category. 
\end{setup}

%%%%%%%%%%%%%%%%%%%%%%%%%%%%%%%%%%%%%%%%%
%%%%%%%%%%%%%%%%%%%%%%%%%%%%%%%%%%%%%%%%%

\subsection{\texorpdfstring{$\BE$}{E}-stratifying systems}
\label{subsec:E-stratifying-systems}

Our first main definition of this section is as follows, and it initiates our study of stratifying systems in extriangulated categories.

\begin{definition}
\label{def:stratifying-system}
Suppose 
$\Phi 
    = \{ \Phi_{i} \}_{i=1}^{n}
    = \{\Phi_{1}, \ldots, \Phi_{n} \}$ 
is a finite collection of indecomposable objects of $\A$.
We call $\Phi$ an \emph{$\BE$-stratifying system (in $(\A,\BE,\fs)$)} if: 
\begin{enumerate}[label=\textup{(S\arabic*)}]
\item\label{S1} $\A(\Phi_{j},\Phi_{i}) =0$  whenever $j > i$; and 
\item\label{S2} $\BE(\Phi_{j},\Phi_{i}) =0$ whenever $j \geq i$. 
\end{enumerate}
\end{definition}

Given an $\BE$-stratifying system $\Phi$ in $\A$, we will study the extriangulated subcategory $\F(\Phi)$ of $\Phi$-filtered objects of $\A$. 
We can immediately prove that 
$\F(\Phi)$ has a `weak' Harder-Narasimhan filtration property 
\cite{%
Harder1975,
Reineke2003,
Chen2010,
Treffinger2018,
tattar2020torsion}. 
For a collection $\Psi$ of indecomposable objects of $\A$, we denote by $\add\Psi$ the full subcategory of $\A$ given by taking the closure of $\Psi$ under direct sums and isomorphisms.

\begin{prop} 
\label{prop:existsnonincreasingfiltration} 
Let 
$\Phi
    = \{ \Phi_{i} \}_{i=1}^{n}
$ 
be a finite collection of indecomposable objects of $\A$ satisfying \ref{S2}. 
The following statements hold. 
\begin{enumerate}[label=\textup{(\roman*)}]

	\item\label{item:existsnonincreasingfiltration-1} For each $1\leq i \leq n$, we have $\F(\Phi_{i}) = \add \Phi_{i}$ is closed under direct summands.

	\item\label{item:existsnonincreasingfiltration-2} Let $X \in \F(\Phi)$. Then for all $\Phi$-filtrations 
$(\xi_{0}, \xi_{1}, \dots, \xi_{t})$ 
\begin{equation}\label{eqn:Phi-filtration}
\begin{tikzcd}[column sep=0.5cm, row sep=0.5cm]
0 \arrow{rr}
&& 0 = X_{0} \arrow{rr}\arrow{dl}
&& X_{1} \arrow{r}\arrow{dl}
%&& X_{2} \arrow{r}\arrow{dl}
& \cdots \arrow{r}
& X_{t-1} \arrow{rr}
&& X_{t} = X \arrow{dl}\\
&0 \arrow[dashed]{ul}{\xi_{0}}
&&\Phi_{j_{1}}\arrow[dashed]{ul}{\xi_{1}}
%&&\Phi_{j_{2}}\arrow[dashed]{ul}{\xi_{2}} 
&&
&&\Phi_{j_{t}} \arrow[dashed]{ul}{\xi_{t}}
&
\end{tikzcd}
\end{equation}
of $X$, 
there exists a $\Phi$-filtration $(\xi'_0, \xi'_1, \dots, \xi'_t)$ of the same length 
with the same $\Phi$-factors 
$\Phi_{j'_1},\ldots,\Phi_{j'_t}$ but re-ordered  
so that 
$j'_{1} \geq j'_{2} \geq \dots \geq j'_{t}$. 
Moreover, 
\begin{equation}\label{eqn:F-Phi-is-star-cat}
\F(\Phi) 
	= \F(\Phi_{n}) \ast \F(\Phi_{n-1}) \ast \dots \ast \F(\Phi_1)
	= \add \Phi_{n} \ast \add \Phi_{n-1} \ast \cdots \ast \add \Phi_{1}.
\end{equation}
\end{enumerate}
\end{prop}

\begin{proof}
\ref{item:existsnonincreasingfiltration-1}:\;\;
It follows from \ref{S2} that $\F(\Phi_{i}) = \add \Phi_{i}$ for each $1\leq i \leq n$, and $\add \Phi_{i}$ is clearly closed under direct summands.

\ref{item:existsnonincreasingfiltration-2}:\;\;
For $t=0,1$, there is nothing to show. 
Thus, suppose $t\geq 2$. 
If $j_{1} \geq j_{2} \geq \dots \geq j_{t}$, then we are done. 
Otherwise, let $s\in\{1,\ldots,t-1\}$ be the largest integer such that $j_{s} < j_{s+1}$.
Applying (ET$4$) to the extriangles $\xi_{s}$ and $\xi_{s+1}$ yields 
\[
\begin{tikzcd}[row sep=0.5cm,row sep=0.6cm]
X_{s-1} \arrow[r] \arrow[d, equal] & X_{s} \arrow[r] \arrow[d]         & \Phi_{j_{s}} \arrow[d]  \arrow[r, dashed, "{\xi_{s}}"] & {}
\\
X_{s-1} \arrow[r]              & X_{s+1} \arrow[r] \arrow[d] & E \arrow[d] \arrow[r, dashed, "\eta"] & {} \\
                             & \Phi_{j_{s+1}} \arrow[r, equal]                       \arrow[d, dashed, "{\xi_{s+1}}"] & \Phi_{j_{s+1}}.     \arrow[d, dashed, "\delta"] \\ &{} & {} & {}
\end{tikzcd}
\]
Since $\Phi$ satisfies \ref{S2}, we have $\BE(\Phi_{j_{s+1}}, \Phi_{j_{s}})=0$, and thus 
$\delta 
	= {}_{\Phi_{j_{s}}}0_{\Phi_{j_{s+1}}}$ 
Without loss of generality, we may assume $E = \Phi_{j_{s+1}} \oplus \Phi_{j_{s}}$. 
We now apply (ET$4)^{\op}$ to $\eta$ and 
%the split extriangle 
${}_{\Phi_{j_{s}}}0_{\Phi_{j_{s+1}}}$ 
to obtain 
\[
\begin{tikzcd}[row sep=0.5cm,row sep=0.6cm]
X_{s-1} \arrow[r] \arrow[d, equal] & X'_{s} \arrow[r] \arrow[d]         & \Phi_{j_{s+1}} \arrow[d]  \arrow[r, dashed, "\xi'_{s}"] & {}
\\
X_{s-1} \arrow[r ]              & X_{s+1} \arrow[r] \arrow[d] & \Phi_{j_{s+1}} \oplus \Phi_{j_{s}} \arrow[d] \arrow[r, dashed, "\eta"] & {} \\
                             & \Phi_{j_{s}} \arrow[r, equal]                       \arrow[d, dashed, "\xi'_{s+1}"] & \Phi_{j_{s}}.      \arrow[d, dashed] \\ &{} & {} & {}               
\end{tikzcd} 
\] 
Now replace $\xi_{s}$ and $\xi_{s+1}$ by $\xi'_s$ and $\xi'_{s+1}$, respectively. 
Note that this yields a $\Phi$-filtration 
$(\xi_{0}, \ldots, \xi_{s-1}, \xi'_{s}, \xi'_{s+1}, \xi_{s+2}, \dots, \xi_{t})$, which 
has exactly the same $\Phi$-factors as $(\xi_{0}, \xi_{1}, \dots, \xi_{t})$ just with $\Phi_{j_{s}}$ and $\Phi_{j_{s}+1}$ switched. 
By repeating this process we obtain the desired $\Phi$-filtration of $X$. 

For the last assertion, note that the $\Phi$-filtration $(\xi'_0, \xi'_1, \dots, \xi'_t)$ implies 
$X\in\add\Phi_{j'_1} \ast \cdots \ast \add\Phi_{j'_t}$ with $j'_{1} \geq j'_{2} \geq \dots \geq j'_{t}$. Therefore, for any $X\in\F(\Phi)$, we have 
$X\in \add \Phi_{n} \ast \dots \ast \add \Phi_1$ since each $\add \Phi_{i}$ is extension-closed by \ref{S2}. 
Equation \eqref{eqn:F-Phi-is-star-cat} then follows from \ref{item:existsnonincreasingfiltration-1}.
\end{proof}

\begin{cor}
\label{cor:closed-under-summands}
Let 
$\Phi
    = \{ \Phi_{i} \}_{i=1}^{n}
$ 
be an $\BE$-stratifying system in  $(\A,\BE,\fs)$. Then $\F(\Phi)$ is closed under direct summands.
\end{cor}
\begin{proof}
By \cref{prop:existsnonincreasingfiltration}, we have 
$\F(\Phi) 
	= \add \Phi_{n} \ast \cdots \ast \add \Phi_{1}
%	= \F(\Phi_{n}) \ast \cdots \ast \F(\Phi_{1})
$ 
and that each $\add(\Phi_{i})$ is closed under direct summands. 
By \ref{S1} 
we have $\A(\add \Phi_{i}, \add \Phi_{j}) = 0$ for all $i > j$.
Thus, $\F(\Phi)$ is closed under direct summands by \cref{lem:summandslemma}.
\end{proof}

We close this subsection by comparing our notion of $\BE$-stratifying system with those in other contexts in the literature already.

\begin{remark} \label{rem-definitioncomparison} 
Let $\Phi 
    = \{ \Phi_{i} \}_{i=1}^{n}$ be a collection of indecomposable objects of $\A$.
    \begin{enumerate}[label=\textup{(\roman*)}]

\item\label{item:triangulated-ESS} Suppose that $(\A,\BE,\fs)$ is a triangulated category with shift functor $[1]$. 
Then $\Phi$ is a \emph{$\Theta$-system} in the sense of \cite[Def.\ 5.1]{MS} if it satisfies the axioms \ref{S1}, \ref{S2} and also
\begin{enumerate}[label=\textup{(S3)}, wide=0pt, leftmargin=35pt, labelwidth=30pt, labelsep=5pt, align=right]
\item\label{S3}
$\A(\Phi_{i}, \Phi_{j} [-1]) = 0 \text{ for all } 1 \leq i,j \leq n.$
\end{enumerate}
We note that although it is not assumed that $\Phi_{i}$ is indecomposable for all $1\leq i \leq n$ in \cite{MS}, if $\Phi$  is a $\Theta$-system then it follows from \cite[Thm.\ 5.9, Prop.\ 6.2]{MS} that each $\Phi_{i}$ is indeed indecomposable. 
Therefore, \cref{def:stratifying-system} is more general than \cite[Def.\ 5.1]{MS}; see Section~\ref{sec:examples}. 
Furthermore, \ref{S3} implies that the subcategory $\F(\Phi)$ is an exact category by \cite[Thm.\ 1]{Dyer-exact-subcategories-of-triangulated-categories}. 
Thus, in some sense, the investigation of $\Theta$-systems in \cite{MS} is a study of exact subcategories of triangulated categories. 
Such subcategories 
have been of interest classically (e.g.\ \cite{BeilinsonBernsteinDeligne-perverse-sheaves}), as well as attracted attention recently 
(e.g.\ \cite{Jorgensen-abelian-subcategories-of-triangulated-categories-induced-by-simple-minded-systems}, \cite[Exam.\ 5.6]{Bennett-TennenhausHauglandSandoyShah-the-category-of-extensions-and-a-characterisation-of-n-exangulated-functors}). 

\item \label{item:exact-ESS} Suppose that $(\A,\BE,\fs)$ is an exact category $(\A,\E)$. Then $\Phi$ is an \emph{$\Ext_{\E}$-system} in the sense of \cite[Def.\ 6.1]{Sa} if and only if the axioms \ref{S1} and \ref{S2} are satisfied.  Therefore, in the exact setting, the notions of an $\BE$-stratifying system and an  $\Ext_{\E}$-system coincide.

\item 
Independently, Adachi--Tsukamoto \cite{adachi2022mixed} have recently investigated similar concepts in the extriangulated setting, namely those of  \emph{mixed standardisable sets} and \emph{mixed (bi)strati\-fying systems}, which generalise work of Dlab--Ringel \cite{DlabRingel-the-module-theoretical-approach-to-quasi-hereditary-algebras}, \cite{ES} and the works discussed above. In \cite{adachi2022mixed}, when working in an extriangulated category $(\A, \BE, \fs)$ equipped with a negative first extension structure $\BE^{-1}$ in the sense of \cite{adachi2021intervals}, 
the authors define and study sets of objects $\Phi = \{\Phi_{1}, \dots, \Phi_{n}\}$ that satisfy similar axioms to those of our $\BE$-stratifying systems and crucially one more condition, namely that $\BE^{-1}(\Phi, \Phi)=0$. This extra assumption implies that $\F(\Phi)$ is an exact category. Thus, the notions we discuss in this paper are more general, but this extra condition is perfectly suitable for the purposes of \cite{adachi2022mixed}.
\end{enumerate}
\end{remark}

%%%%%%%%%%%%%%%%%%%%%%%%%%%%%%%%%%%%%%
%%%%%%%%%%%%%%%%%%%%%%%%%%%%%%%%%%%%%%

\subsection{Projective \texorpdfstring{$\BE$}{E}-stratifying systems}
\label{subsec:weak-projective-stratifying-systems}

In this subsection we discuss projective objects relative to an $\BE$-stratifying system. 
Our main result here is that each $\BE$-stratifying system can be extended to a projective $\BE$-stratifying system (see \cref{def:weak-projective-system} and \cref{thm:existence-of-proj-ESS}). 
We still assume \Cref{setup:A-is-artin-extriangulated}.

\begin{notation}
\label{notation:Theta-subscript-sets}
If $\Phi = \{\Phi_{i}\}_{i=1}^{n}$ is an $\BE$-stratifying system and $i_{0}\in\{1,\ldots,n\}$, then we denote the subset 
$\set{\Phi_{j} | j >i_{0}}$ by $\Phi_{j>i_{0}}$. 
We define other subsets (e.g.\ 
$
\Phi_{j\leq i_{0}} 
    \deff 
    \set{ \Phi_{j} | j \leq i_{0} } 
$) 
of $\Phi$ in a similar way.
\end{notation}

Our second key definition of \cref{sec:stratifying-systems} is the following.

\begin{definition}
\label{def:weak-projective-system} 
A \emph{projective $\BE$-stratifying system (in $(\A,\BE,\fs)$)} is a pair $(\Phi,Q)$, where $\Phi = \{ \Phi_{i} \}_{i=1}^{n}$ is an $\BE$-stratifying system and 
$
Q 
    = \{ Q_{i} \}_{i=1}^{n} 
    = \{Q_{1},\ldots, Q_{n}\}
$ 
is a set of 
objects of $\A$, satisfying the following conditions. 
\begin{enumerate}[label=\textup{(PS\arabic*)}]

\item\label{PS1} $\BE(Q_{i},\Phi_{j}) = 0$ for all $1\leq i,j\leq n$. 

\item\label{PS2} For each $1\leq i \leq n$, there exists an 
extriangle 
\stepcounter{equation}
\begin{equation}
\label{eqn:PS2-extriangles}
\begin{tikzcd}%[column sep=0.5cm]
K_{i} \arrow{r}{k_{i}} & Q_{i} \arrow{r}{q_{i}} & \Phi_{i} \arrow[dashed]{r}{\eta_{i}}& {},
\end{tikzcd}
\tag*{$(\theequation)_{i}$}
\end{equation}
in $(\A,\BE,\fs)$ with 
$K_{i}$ in $\F(\Phi_{j>i})$. 
\end{enumerate}
We say that a projective $\BE$-stratifying system 
$(\Phi,Q)
    = (\{\Phi_{i}\}_{i=1}^{n},\{Q_{i}\}_{i=1}^{n})$
is \emph{minimal} if, for each $1\leq i \leq n$, 
the morphism $q_{i}$ in the extriangle \ref{eqn:PS2-extriangles} is right minimal.
\end{definition}

\begin{remark}
\label{rem:on-def-projective-ESS}
Let $(\Phi,Q)
    = (\{\Phi_{i}\}_{i=1}^{n},\{Q_{i}\}_{i=1}^{n})$ 
be a projective $\BE$-stratifying system. 
\begin{enumerate}[\textup{(\roman*)}]

    \item\label{item:kn-is-zero} 
    Since $\Phi_{j>n} = \emptyset$ and $\F(\emptyset) = \{0\}$ (see \cref{rem:comments-on-FTheta}\ref{part:F-empty-is-zero}), we have that $K_{n} = 0$ and so 
    $q_{n}\colon Q_{n} \overset{\iso}{\longrightarrow} \Phi_{n}$ is a non-zero isomorphism as $\Phi_{n} \neq 0$. 
       
    \item In contrast to the triangulated \cite{MS} and exact \cite{Sa} cases, we allow the possibility that $Q_{i} \cong 0$ for $i<n$; 
    see \cref{example:zero projective}.

    \item Since $K_{i}$ and $\Phi_{i}$ lie in the extension-closed subcategory $\F(\Phi_{j\geq i})\sse \F(\Phi)$ (see 
    \cref{rem:comments-on-FTheta}\ref{part:F(X)-is-smallest-ext-closed-subcat-containing-X}), we have  $Q_{i}\in\F(\Phi_{j\geq i})\sse \F(\Phi)$ also. Thus, each extriangle \ref{eqn:PS2-extriangles} from \ref{PS2} is an extriangle in $(\F(\Phi),\restr{\BE}{\F(\Phi)},\restr{\fs}{\F(\Phi)})$.
    
    \item\label{part:Q_{i}-is-projective}
    An \emph{$\BE'$-projective} object $P$ in an extriangulated category $(\A',\BE',\fs')$ is an object that satisfies 
    $\BE'(P,A)=0$ for all $A\in\A'$; see \cite[Def.\ 3.23, Prop.\ 3.24]{NP19}. 
    We denote the full subcategory of $\A$ consisting of all objects that are projective relative to $\F(\Phi)$ as  
    \[
    \cP(\Phi) 
        \deff 
        \set{ Z \in \A | \BE(Z, X) = 0  \text{ for all  } X \in \F(\Phi)}. 
    \]
    Thus, the class of objects $\cP(\Phi) \cap \F(\Phi)$ is the collection of all $\restr{\BE}{\F(\Phi)}$-projective objects in $(\F(\Phi),\restr{\BE}{\F(\Phi)},\restr{\fs}{\F(\Phi)})$.

    It is easy to verify that $q_{i}\colon Q_{i} \to \Phi_{i}$ is a right $\cP(\Phi)$-approximation of $\Phi_{i}$ since $K_{i}\in\F(\Phi)$
    Condition \ref{PS1} tells us that each $Q_{i}$ is $\restr{\BE}{\F(\Phi)}$-projective in $(\F(\Phi),\restr{\BE}{\F(\Phi)},\restr{\fs}{\F(\Phi)})$ by \cref{lem:Zhou-lem3-3-orthogonality-compatible-with-F}\ref{part:Zhou3.3ii}, i.e.\ $Q_{i}\in\cP(\Phi)\cap\F(\Phi)$.

\end{enumerate}
\end{remark}

Continuing from \cref{rem-definitioncomparison}, let us compare \cref{def:weak-projective-system} with the versions in the literature for triangulated and exact categories.

\begin{remark} 
\label{rem-definitioncomparison2}
Let $(\Phi,Q)
    = (\{\Phi_{i}\}_{i=1}^{n},\{Q_{i}\}_{i=1}^{n})$ 
be a pair of sets of indecomposable objects of $\A$.
\begin{enumerate}[\textup{(\roman*)}]
    \item Suppose that $(\A,\BE,\fs)$ is a triangulated category with shift functor $[1]$. Then $(\Phi, Q)$ is a  \emph{projective $\Theta$-system} in the sense of \cite[Def.\ 5.2]{MS}  if  $\Phi$ is a $\Theta$-system and axioms equivalent to \ref{PS1} and \ref{PS2} are satisfied. 
    It follows from the discussion in \cref{rem-definitioncomparison}\ref{item:triangulated-ESS} that our projective $\BE$-stratifying systems do indeed generalise projective $\Theta$-systems in this setting. 

    \item Suppose that $(\A,\BE,\fs)$ is an exact category $(\A,\E)$. Then $(\Phi, Q)$ is a \emph{projective $\Ext_{\E}$-system} in the sense of \cite[Def.\ 6.2]{Sa} if and only if $\Phi$ is an $\BE$-stratifying system  and the axioms \ref{PS1} and \ref{PS2} are satisfied. Hence, for exact categories, the notions of  projective $\BE$-stratifying system and projective $\Ext_{\E}$-system  coincide.  
\end{enumerate}
\end{remark}

The next lemma provides an analogue for an extriangulated category of the `universal extension' in a module category \cite[Lem.\ VIII.2.4]{Skoronski-Yamagata-II}. 
It is due to \cite[Lem.\ 2.12]{Zh19}; note, however, that there is a typo in that statement. 
The proof of \cite[Lem.\ 3.4]{MS} 
readily translates over. 
However, we remark that the additivity of the category of $\BE$-extensions $\BE$-$\Ext(\A)$ (see \cref{sec:extriangulated-categories}) is needed for this proof.

Recall that $\ring$ is a commutative artinian ring. We denote the \emph{length} of an $\ring$-module $X$ by $\ell_{\ring}(X)$.

\begin{lem}
\label{lem:Zhou-lemma-2-12}
Assume $A,C\in\A$ with $0 < l\deff \ell_{\ring}(\BE(C,A)) <\infty $. Then there exists an extriangle 
$\begin{tikzcd}[column sep=0.5cm,cramped]
A^{\oplus l} \arrow{r}& B \arrow{r}& C \arrow[dashed]{r}{\eta}& {},
\end{tikzcd}$ 
such that $\eta\neq 0$ and $\eta_{A}^{\sharp}\colon \A(A^{\oplus l},A) \to \BE(C,A)$ is surjective. Furthermore, if $\BE(A,A)=0$, then $\BE(B,A)=0$.
\end{lem}

Recall that an extriangulated category $(\A',\BE',\fs')$ \emph{has enough $\BE'$-projectives} if for each $X\in\A'$ there is an extriangle 
$
\begin{tikzcd}[column sep=0.5cm,cramped]
A \arrow{r}& P \arrow{r}& X \arrow[dashed]{r}&{} 
\end{tikzcd}
$
in $(\A',\BE',\fs')$ where $P$ is $\BE'$-projective (see \cite[Def.\ 3.25]{NP19}). 
Compare the next result with \cite[Lems.\ 4.12, 4.13]{MS}.

\begin{prop}
\label{thm: ext proj approx}
Suppose $\Phi = \{\Phi_{i}\}_{i=1}^{n}$ is an $\BE$-stratifying system. 
Then for each $X \in \A$, there exists an extriangle
$
\begin{tikzcd}[column sep=0.5cm,cramped]
A \arrow[r] & B \arrow[r] & X \arrow[r, dashed, "\eta"] & {} 
\end{tikzcd}
$
with $A \in \F(\Phi)$ and $B \in \cP(\Phi)$.
In particular, the extriangulated category $(\F(\Phi),\restr{\BE}{\F(\Phi)},\restr{\fs}{\F(\Phi)})$ has enough $\restr{\BE}{\F(\Phi)}$-projectives.
\end{prop}

\begin{proof}
We will construct extriangles 
$\begin{tikzcd}[cramped]
\Phi_{a_{j}}^{\oplus m_{j}} \arrow[r] & X_{j+1} \arrow[r, "{y_{j+1}}"] & X_{j} \arrow[dashed, "{\delta_{j}}"]{r} & {}
\end{tikzcd}$
with $0 \leq j \leq r$ for some $r\in\{0,\ldots,n\}$ dependent on $X$, 
such that $m_{j} \geq 0$, and:
\begin{enumerate}[label=\textup{(\arabic*)}]
    \item $1\leq a_{j} \leq n$ is minimal such that $\BE(X_{j},\Phi_{a_{j}}^{\oplus m_{j}}) \neq 0$ and $\BE(X_{j},\Phi_{s}) = 0$ for all $s<a_{j}$;
    \item $a_{j+1}>a_{j}$ for $0\leq j \leq r-1$; and 
    \item $\BE(X_{r+1}, \Phi) = 0$.
\end{enumerate}

Set $X_{0}\deff X$. If $\BE(X_{0},\Phi) = 0$, then take $\delta_{0}$ to be the extriangle 
$\begin{tikzcd}[column sep=0.9cm,cramped]
0 \arrow{r}& X_{0} \arrow{r}{\id{X_{0}}}& X_{0} \arrow[dashed]{r}& {}
\end{tikzcd}$
and stop with $X_{1} \deff X_{0} = X$. 
Else, take $a_{0} \in \{1,\ldots, n\}$ to be minimal such that $\BE(X_{0},\Phi_{a_{0}}) \neq 0$. Apply \cref{lem:Zhou-lemma-2-12} to obtain a non-split extriangle
$\begin{tikzcd}[column sep=0.6cm,cramped]
\Phi_{a_{0}}^{\oplus m_{0}} \arrow{r}& X_{1} \arrow{r}{y_{1}}& X_{0} \arrow[dashed]{r}{\delta_{0}} &{}
\end{tikzcd}$ 
such that $(\delta_{0}^{\sharp})_{\Phi_{a_{0}}} \colon  \A(\Phi_{a_{0}}^{\oplus m_{0}}, \Phi_{a_{0}}) \to \BE(X, \Phi_{a_{0}})$ is surjective. 

If $\BE(X_{1},\Phi) = 0$, then we are done with $r=0$.
Else, 
set 
\[
a_{1}
    \deff 
    \text{min}\Set{ 1\leq s \leq n | 
        \BE(X_{1},\Phi_{s}) \neq 0\text{ and }
        \BE(X_{1}, \Phi_{s'}) =0 \;\; \forall s' <s  }.
\] 
There is an extriangle 
$\begin{tikzcd}[column sep=0.5cm,cramped]
\Phi_{a_{1}}^{\oplus m_{1}} \arrow{r} & X_{2} \arrow{r}{y_{2}} & X_{1} 
\arrow[dashed, "{\delta_{1}}"]{r} & {}
\end{tikzcd}$ 
such that $(\delta_{1}^{\sharp})_{\Phi_{a_{1}}} \colon  \A(\Phi_{a_{1}}^{\oplus m_{1}}, \Phi_{a_{1}}) \to \BE(X_{1}, \Phi_{a_{1}})$ is surjective by \cref{lem:Zhou-lemma-2-12}.

Assume for contradiction that $a_{1} \leq a_{0}$. 
The extriangle $\delta_{0}$ induces an exact sequence
\[
\begin{tikzcd}[column sep=1.7cm]
\A(\Phi_{a_{0}}^{\oplus m_{0}}, \Phi_{a_{1}}) 
    \arrow{r}{(\delta_{0}^{\sharp})_{\Phi_{a_{1}}}}
&\BE(X_{0}, \Phi_{a_{1}}) 
    \arrow{r}{\BE(y_{1}, \Phi_{a_{1}})}
&\BE(X_{1}, \Phi_{a_{1}}) 
    \arrow{r}
&\BE(\Phi_{a_{0}}^{\oplus m_{0}}, \Phi_{a_{1}}).
\end{tikzcd} 
\]
If $a_{1} < a_{0}$, then  
$\BE(X_{0}, \Phi_{a_{1}}) = 0$ by the minimality of $a_{0}$ 
and 
$\BE(\Phi_{a_{0}}^{\oplus m_{0}}, \Phi_{a_{1}}) = 0$ by \ref{S2}.  
If $a_{1} =a_{0}$, then we would still have 
$\BE(\Phi_{a_{0}}^{\oplus m_{0}}, \Phi_{a_{1}}) = 0$ by \ref{S2}, and now also that 
$\BE(y_{1}, \Phi_{a_{1}}) = 0$ since $(\delta_{0}^{\sharp})_{\Phi_{a_{1}}} = (\delta_{0}^{\sharp})_{\Phi_{a_{0}}}$ is surjective. 
In either case, we deduce that $\BE(X_{1}, \Phi_{a_{1}}) = 0$, which is a contradiction, so we must have $a_{1} > a_{0}$.

Then we can repeatedly apply the above argument to construct the desired extriangles, noting that this process will terminate since $\{1,\ldots,n\}$ is finite. 
Lastly, by repeatedly applying (ET4$)^{\op}$ to the deflations $y_{j}$, we obtain an extriangle 
$
\begin{tikzcd}[column sep=0.5cm,cramped]
A \arrow{r} & X_{r+1} \arrow{r}{y} & X \arrow[dashed]{r}{\eta} & {}
\end{tikzcd}
$
with $y\deff y_{1}y_{2}\dots y_{r+1}$, 
$A 
    \in \F(\Phi_{j\geq a_{0}})
    \sse \F(\Phi)
$ 
and 
$\BE(X_{r+1},\Phi) = 0$.

The last assertion follows from $\F(\Phi)$ being  extension-closed in $(\A,\BE,\fs)$, noting that if $X\in\F(\Phi)$ then $X_{r+1}$ must also lie in $\F(\Phi)$.
\end{proof}

We are now in position to prove that each $\BE$-stratifying system is always part of a minimal projective one. This is the main result of this subsection.

\begin{thm}
\label{thm:existence-of-proj-ESS}
Let $\Phi= \{\Phi_{i}\}_{i=1}^{n}$ be an $\BE$-stratifying system. 
Then there exists a set $Q = \{Q_{i}\}_{i=1}^{n}$ of objects of $\A$, unique up to isomorphism, such that $(\Phi, Q)$ is a minimal projective $\BE$-stratifying system. 
\end{thm}

\begin{proof} 
We first show that each $\BE$-stratifying system $\Phi= \{\Phi_{i}\}_{i=1}^{n}$ is part of a projective one. Then we will show that we may refine this projective $\BE$-stratifying system to a minimal one.
Apply \cref{thm: ext proj approx} to $\Phi_{1},\Phi_{2},\ldots,\Phi_{n}$ in that order to obtain 
a sequence of objects $Q_{1},Q_{2},\ldots,Q_{n}$. 
Thus, for each $1\leq i \leq n$, we have an extriangle 
\begin{equation}\label{eqn:Phi-porj-system-extriangle}
\begin{tikzcd}%[column sep=0.5cm]
K_{i}\arrow{r} & Q_{i} \arrow{r}{q_{i}}& \Phi_{i} \arrow[dashed]{r}{\eta_{i}}&{}
\end{tikzcd}
\end{equation}
with $Q_{i}$ projective in $(\F(\Phi),\restr{\BE}{\F(\Phi)},\restr{\fs}{\F(\Phi)})$. 
Furthermore, the fact that $K_{i}$ lies in $\F(\Phi_{j>i})$ is 
a consequence of the construction of $a_{0}$ in the proof of \cref{thm: ext proj approx}, noting that \ref{S2} implies $i < a_{0}$.
Moreover, with $Q \deff \{Q_{i}\}_{i=1}^{n}$, the pair 
$(\Phi, Q)$ is a projective $\BE$-stratifying system.

We now show that, for each $i \in \{1, \dots, n\}$, we may take the morphism $q_{i}$ in \eqref{eqn:Phi-porj-system-extriangle} to be right minimal.
If $q_{i} = 0$ is zero, then by \cref{lem:decomposable} we have that 
\[
 \eta_i 
    \cong 
    \langle 
        \begin{tikzcd}
        K_{i}' \arrow{r} & 0 \arrow{r}{} & \Phi_{i},
        \end{tikzcd}
    \eta_i' 
    \rangle 
    \oplus 
    \langle
        \begin{tikzcd}
        Q_{i} \arrow{r}{\id{Q_{i}}} & Q_{i} \arrow{r} & 0,
        \end{tikzcd} 
    {}_{Q_{i}}0_0 
    \rangle, 
\] 
where $K'_{i} \oplus Q_{i}\iso K_{i}$. 
Note that the morphism $q'_{i} \colon 0\to \Phi_{i}$ is a monomorphism and hence right minimal.
Thus, in this case, we can replace $\eta_i$ by $\eta_i'$, noting that $K_{i}' \in \F(\Phi_{j>i})$ 
as $\F(\Phi_{j>i})$ is closed under direct summands 
by \cref{cor:closed-under-summands}.

Now suppose $q_{i}$ is non-zero. 
Since $\A$ is a Krull-Schmidt category and $q_{i}$ is a right $\cP(\Phi)$-approximation of $\Phi_{i}$ (see \cref{rem:on-def-projective-ESS}\ref{part:Q_{i}-is-projective}), there is a direct sum decomposition $Q_{i} \cong Q'_{i} \oplus P$ such that $q'_{i}\deff \restr{q_{i}}{Q'_{i}} \neq 0$ is right minimal and $\restr{q_{i}}{P} = 0$ by \cite[Cor.~1.4]{KrauseSaorin-minimal-approximations-of-modules}. 
Again by \cref{lem:decomposable}, we have  
\[
 \eta_i \cong \langle \begin{tikzcd}
 K_{i}' \arrow{r} & Q'_{i} \arrow{r}{q'_{i}} & \Phi_{i},
 \end{tikzcd} \eta_i' \rangle \oplus \langle  \begin{tikzcd}
 P \arrow{r}{\id{P}} & P \arrow{r} & 0,
 \end{tikzcd} {}_{P}0_0 \rangle,
\] and we replace $\eta_i$ by $\eta_i'$ as before.

Thus, we have obtained a minimal $\BE$-projective system $(\Phi,Q') = (\{\Phi_{i}\}_{i=1}^{n},\{Q'_{i}\}_{i=1}^{n})$. 
Since we now have that $q'_{i}\colon Q'_{i} \to \Phi_{i}$ is a minimal right $\cP(\Phi)$-approximation of $\Phi_{i}$ (see \cref{rem:on-def-projective-ESS}\ref{part:Q_{i}-is-projective}), we deduce that $Q'_{i}$ is uniquely determined up to isomorphism.
\end{proof}

%%%%%%%%%%%%%%%%%%%%%%%%%%%%%%%%%%%%%%%%%%%%%%%%%%%%%%%
%%%%%%%%%%%%%%%%%%%%%%%%%%%%%%%%%%%%%%%%%%%%%%%%%%%%%%%

\subsection{The Jordan-H{\"{o}}lder property and projective \texorpdfstring{$\BE$}{E}-stratifying systems}

In this subsection we identify a sufficient condition for a minimal projective $\BE$-stratifying system to give rise to a Jordan-H{\"{o}}lder extriangulated category (see \cref{def:E-s-notions}\ref{item:Jordan-Holder-extrianguated-category}). We need two preliminary lemmas.

\begin{lem}
\label{lem:non-homs-Qi-to-Thetai}
Let $(\Phi,Q) = (\{\Phi_{i}\}_{i=1}^{n},\{Q_{i}\}_{i=1}^{n})$ be a  projective $\BE$-stratifying system. 
If $i>j$, then $\A(Q_{i},\Phi_{j})=0$.
\end{lem}
\begin{proof}
The extriangle 
\ref{eqn:PS2-extriangles}
induces an exact sequence
$
\begin{tikzcd}[column sep=0.5cm,cramped]
\A(\Phi_{i},\Phi_{j}) \arrow{r}& \A(Q_{i},\Phi_{j})\arrow{r} & \A(K_{i},\Phi_{j}).
\end{tikzcd}
$ 
If $i>j$, then 
$\A(\Phi_{i},\Phi_{j})=0$ by \ref{S1}. 
Furthermore, $\A(K_{i},\Phi_{j})=0$ using that  $K_{i}\in\F(\Phi_{l>i})$, where $l>i>j$, and \ref{S1} and \cref{lem:Zhou-lem3-3-orthogonality-compatible-with-F}\ref{part:Zhou3.3i}. 
Hence, $\A(Q_{i},\Phi_{j})=0$.
\end{proof}

For the next lemma, we need some terminology first. 

\begin{definition}
\label{def:exact-functor}
Suppose $(\A',\BE',\fs')$ is an extriangulated category and that $F\colon \A' \to \Ab$ is an additive, covariant functor. 
We call $F$ \emph{left exact on $(\A',\BE',\fs')$} if, for each extriangle 
$
\begin{tikzcd}[column sep=0.5cm,cramped]
A \arrow{r}{a}& B \arrow{r}{b}& C \arrow[dashed]{r}& {}
\end{tikzcd}
$
in $\A$, 
the sequence
$\begin{tikzcd}[column sep=0.7cm,cramped]
0 \arrow{r}& FA \arrow{r}{Fa}& FB
\end{tikzcd}$
is exact in $\Ab$. 
Being \emph{right exact on $(\A',\BE',\fs')$} is defined dually. 
Lastly, we say that $F$ is \emph{exact on $(\A',\BE',\fs')$} if $F$ is both left and right exact on $(\A',\BE',\fs')$. 
\end{definition}

Note here that if we equip $\Ab$ with its canonical extriangulated structure given that it is an abelian category, then an additive, covariant functor $F\colon \A' \to \Ab$ is exact on $(\A',\BE',\fs')$ in the sense above if and only if it is part of an \emph{extriangulated} functor in the sense of \cite[Def.\ 2.32]{Bennett-TennenhausShah-transport-of-structure-in-higher-homological-algebra}.

\begin{lem} \label{lem:left-exact-no-zero-proj}
Let $(\Phi, Q) = (\{\Phi_{i}\}_{i=1}^{n},\{Q_{i}\}_{i=1}^{n})$ be a minimal projective $\BE$-stratifying system. If $\A(Q_{i},-)$ is left exact on $(\F(\Phi),\restr{\BE}{\F(\Phi)},\restr{\fs}{\F(\Phi)})$ for each $1 \leq i \leq n $, then in the extriangles \ref{eqn:PS2-extriangles} each $q_{i}\colon Q_{i}\to \Phi_{i}$ is non-zero.
\end{lem}
\begin{proof}
By Remark~\ref{rem:on-def-projective-ESS}\ref{item:kn-is-zero} we have that $q_{n}\colon Q_n \to \Phi_{n}$ is a non-zero isomorphism. 
We now proceed by induction. Thus, fix $i<n$ and suppose that $q_j$ is non-zero for all $j\in\{ i+1,\ldots, n \}$. 
Assume for contradiction that $q_{i}\colon Q_{i}\to \Phi_{i}$ is zero, which implies $Q_{i} = 0$ as $q_{i}$ is right minimal by assumption on  $(\Phi, Q)$.
Since $K_{i} \in \F(\Phi_{j>i})$, by \cref{prop:existsnonincreasingfiltration}, we know there is a $\Phi_{j>i}$-filtration $(\xi_{0}, \xi_{1}, \dots, \xi_{t})$ of $X = K_{i}$ as in \eqref{eqn:Phi-filtration}, with and $n\geq j_{1}\geq j_{2}\geq\ldots,\geq j_{t} >i$. 

Let $k\in\{ i+1,\ldots, n \}$ be arbitrary. 
We shall show that $\A(Q_k, X_{l}) = 0$ for all $0\leq l \leq t$. 
Applying $\A(Q_k, -)$ to the extriangle \ref{eqn:PS2-extriangles}, we see that $\A(Q_k, X_{t}) = \A(Q_k, K_{i}) = 0$ as $Q_{i} = 0$. 
Applying $\A(Q_k, -)$ to $\xi_{t}$ we conclude that 
$\A(Q_k, X_{t-1}) =0$ also. 
Continuing in this way we see that $\A(Q_k, \Phi_{j_{1}}) = \A(Q_k, X_{1}) = 0$. 

However, this contradicts that $q_{j_{1}}$ is non-zero. Therefore, we must have that $t=0$ and hence $K_{i} = 0$. From this we see that $\Phi_{i} = 0$, which is our desired contradiction, so that $q_{i}\neq0$.
\end{proof}

We are now in place to present our main result for this subsection, the proof of which 
follows the classical triangular matrix argument that is used in \cite[Lem.\ 1.4]{ES}, \cite[Prop.\ 5.11]{MS} and \cite[Prop.\ 6.5]{Sa}.

\begin{thm}
\label{thm:JHP-for-Theta-filtrations}
Let $(\Phi, Q) = (\{\Phi_{i}\}_{i=1}^{n},\{Q_{i}\}_{i=1}^{n})$ be a minimal projective $\BE$-stratifying system and 
%in the extriangulated category $(\A,\BE,\fs)$. 
suppose that 
$\A(Q_{i},-)$ is left exact on $(\F(\Phi),\restr{\BE}{\F(\Phi)},\restr{\fs}{\F(\Phi)})$ for each $1 \leq i \leq n $. 
Then 
the Jordan-H{\"{o}}lder property holds for $\Phi$-filtrations, that is, any two $\Phi$-filtrations of an object 
$M\in\A$
have the same length and the same $\Phi$-factors, up to isomorphism and permutation.  
\end{thm}

\begin{proof}
 Let $M\in\F(\Phi)$ be arbitrary, and 
let $\xi = (\xi_{0},\ldots,\xi_{t})$
be a $\Phi$-filtration of $M$ as in \eqref{eqn:X-filtration}. 
We apply the functor $\A(Q_{i}, -)$
to each extriangle $\xi_{l}$, for $l = 0,\ldots, t$, which gives a short exact sequence
\[
\begin{tikzcd}
0\arrow{r} & \A(Q_{i},M_{l-1}) \arrow{r}& \A(Q_{i},M_{l}) \arrow{r}& \A(Q_{i},\Phi_{j_{l}}) \arrow{r}&  0,
\end{tikzcd}
\]
using also that $Q_{i} \in \cP(\Phi)$ and that $\A(Q_{i},-)$ is left exact on $\F(\Phi)$ by assumption. 
Recall that under our assumptions, each $\Hom$-set in $\A$ is a finite length $\ring$-module. 
Thus, taking lengths as $\ring$-modules yields an equation 
\stepcounter{equation}
\begin{equation}
\label{equation:exact}
\ell_{\ring}(\A(Q_{i},M_{l}) )
    = \ell_{\ring}(\A(Q_{i},M_{l-1}) )+ \ell_{\ring}(\A(Q_{i},\Phi_{j_{l}}))
\tag*{$(\theequation)_{l}$}
\end{equation}
for each $1\leq l \leq t$. 

Denote by $D = (d_{ij})$ the $n \times n$-matrix with entries
$
d_{ij} \deff \ell_{\ring}(\A(Q_{i},\Phi_{j})).
$ 
By 
\cref{lem:left-exact-no-zero-proj}, 
we have that $q_{i} \neq 0$ and so 
$d_{ii}\neq 0$ for each $1\leq i \leq n$. 
Moreover, $d_{ij} = 0$ for $i>j$ 
by 
\cref{lem:non-homs-Qi-to-Thetai}. 
Thus, $D$ is an upper triangular matrix with non-zero diagonal entries and, hence, $D$ is invertible. 

Further, define  
$c_{i}(X) \deff \ell_{\ring}(\A(Q_{i},X))$ for $X\in\F(\Phi)$, 
and denote by $m_{j} \deff [M : \Phi_{j}]_{\xi}$ the multiplicity of $\Phi_{j}$ in the $\Phi$-filtration $\xi$ of $M$, that is, the number of indices $l$ where $j_{l}=j$.
We claim that the following equation holds for each $i\in\{1,\ldots,n\}$:
\stepcounter{equation}
\begin{equation}
\label{eqn:c_i}
c_{i}(M) = \sum_{j=1}^{n} d_{ij} m_{j}.
\tag*{$(\theequation)_{i}$}
\end{equation}
Indeed, from the equations \ref{equation:exact}, we obtain 
\begin{align*}
c_{i}(M)
    &= c_{i}(M_{t}) \\
    &= \ell_{\ring}(\A(Q_{i},M_{t-1}) ) + \ell_{\ring}(\A(Q_{i},\Phi_{j_{t}}))\\
    &= c_{i}(M_{t-1}) + d_{ij_{t}} \\
    & \hspace{6pt} \vdots \\
    &= d_{ij_{1}} + \cdots + d_{ij_{t}},
\end{align*}
whence \ref{eqn:c_i} follows. 

With 
$
\ul{c} 
    \deff 
     ( c_{1}(M),  c_{2}(M) , \cdots , c_{n}(M))^{T}
$
and 
$
\ul{m}
    \deff 
		( m_{1}, m_{2},\ldots, m_{n})^{T}
$, 
the equations \ref{eqn:c_i} yield the matrix equation
$
D  \ul{m} = \ul{c}.
$
Since $D$ is invertible, we obtain
$ 
\ul{m} = D^{-1} \ul{c}.
$
The right-hand side is independent of the chosen $\Phi$-filtration $\xi$ of $M$, and thus the multiplicities $m_{j}$ of the $\Phi$-factors $\Phi_{j}$ depend only on $M$. 
In particular, this implies that any 
$\Phi$-filtration of $M$ will have length 
$\sum_{j=1}^{n}m_{j}$, 
with $\Phi_{j}$ appearing precisely $m_{j}$ times.
\end{proof}

\begin{cor}
\label{cor:JHP-for-F-Theta}
Let $(\Phi, Q) = (\{\Phi_{i}\}_{i=1}^{n},\{Q_{i}\}_{i=1}^{n})$ be a minimal projective $\BE$-stratifying system and 
suppose that 
$\A(Q_{i},-)$ is left exact on $(\F(\Phi),\restr{\BE}{\F(\Phi)},\restr{\fs}{\F(\Phi)})$ for each $1 \leq i \leq n $. 
Then, up to isomorphism, the collection of $(\restr{\BE}{\F(\Phi)},\restr{\fs}{\F(\Phi)})$-simple objects is precisely the set $\Phi$, and hence 
$(\F(\Phi),\restr{\BE}{\F(\Phi)},\restr{\fs}{\F(\Phi)})$ is a length and Jordan-H{\"{o}}lder extriangulated category.
\end{cor}

\begin{proof}
By \cref{cor:simple-isomorphic-to-X}, we know any $(\restr{\BE}{\F(\Phi)},\restr{\fs}{\F(\Phi)})$-simple object is isomorphic to an object in $\Phi$. For the converse, we will use \cref{lem:0-simple-like-implies-stronger-simple}. Thus, fix $i\in\set{1,\ldots,n}$ and assume there is an extriangle 
$\begin{tikzcd}[column sep=0.5cm,cramped]
    A \arrow{r}{}& \Phi_{i} \arrow{r}& C \arrow[dashed]{r}& {.}
\end{tikzcd}$
We must show $A=0$ or $C=0$. 
Choose a $\Phi$-filtration of $A$ (resp.\ $C$) of length say $u$ (resp.\ $v$). 
By \cref{lem:comp-series-of-A-and-C-gives-one-for-B}, 
we can construct a $\Phi$-filtration of $\Phi_{i}$ of length $u+v$. 
On the other hand, $\Phi_{i}$ has the $\Phi$-filtration 
$( {}_{0}0_{0}, {}_{0}0_{\Phi_{i}} )$
of length $1$. 
But, by \cref{thm:JHP-for-Theta-filtrations}, 
we know the Jordan-H{\"{o}}lder property holds for $\Phi$-filtrations, so 
we must have $u + v = 1$. 
Thus, either $u=0$ and $A=0$, or $v=0$ and $C=0$. 
Hence, $\Phi_{i}$ is an $(\restr{\BE}{\F(\Phi)},\restr{\fs}{\F(\Phi)})$-simple object. 

Each object in $\F(\Phi)$ has a $\Phi$-filtration by definition. 
We have just shown it is an $(\restr{\BE}{\F(\Phi)},\restr{\fs}{\F(\Phi)})$-composition series, so $(\F(\Phi),\restr{\BE}{\F(\Phi)},\restr{\fs}{\F(\Phi)})$ is length. 
It is Jordan-H{\"{o}}lder by \cref{thm:JHP-for-Theta-filtrations}.
\end{proof}

\begin{remark}
We give some examples of situations where the conditions of \cref{thm:JHP-for-Theta-filtrations} hold. 
\begin{enumerate}[label=\textup{(\roman*)}]

\item Clearly the condition that $\A(X, -)$ is left exact on $\F(\Phi)$ for all $X \in \F(\Phi)$ is sufficient for \cref{thm:JHP-for-Theta-filtrations}. Therefore, in light 
of Remarks \ref{rem-definitioncomparison} and \ref{rem-definitioncomparison2}, 
a projective $\Theta$-system in a triangulated category and an $\Ext_{\E}$-system in an exact category give rise to an exact category $\F(\Phi)$ with the Jordan-H{\"{o}}lder property, when the underlying category satisfies \Cref{setup:A-is-artin-extriangulated}.

\item If $\A$ is equipped with a negative first extension structure $\BE^{-1}$ in the sense of \cite{adachi2021intervals}, then the condition of $\A(P_{i}, -)$ being left exact on $\F(\Phi)$ is equivalent to $\BE^{-1}(P_{i}, -)$ being right exact on $\F(\Phi)$. Alternatively,  $\restr{\BE^{-1}(P_{i}, -)}{\F(\Phi)} = 0$ gives a sufficient condition for the left exactness of $\restr{\A(P_{i},-)}{\F(\Phi)}$. 
\end{enumerate}
\end{remark}

Lastly, we comment on the dual of the theory introduced in this section.

\begin{remark}
\label{remark:self-dual}
By dualising the concepts and results of this section, we have that each $\BE$-stratifying system is part of an injective $\BE$-stratifying system which, under dual assumptions, the Jordan-H{\"{o}}lder property holds for $\Phi$-filtrations. 
The definition of an extriangulated category is self-dual \cite[Caution 2.20]{NP19}. 
In other words, if $(\A,\BE,\fs)$ is an extriangulated category, then $(\A^{\op},\BE^{\op},\fs^{\op})$ is extriangulated where $\BE^{\op}$ and $\fs^{\op}$ are defined in the obvious way. In light of this, one may define an
\emph{injective $\BE$-stratifying system in $(\A,\BE,\fs)$} as a pair $(\Phi,Q)
    = (\{\Phi_{i}\}_{i=1}^{n},\{Q_{i}\}_{i=1}^{n})$ of sets of indecomposable objects of $\A$ such that $(\Phi^{\op}, Q^{\op})$ is a projective $\BE^{\op}$-stratifying system in $(\A^{\op},\BE^{\op},\fs^{\op})$. Note that in $(\A^{\op},\BE^{\op},\fs^{\op})$ the axioms of a projective $\BE^{\op}$-stratifying system are considered with respect to the usual ordering of $\mathbb{N}^{\op}$.
\end{remark}

% %%%%%%%%%%%%%%%%%%%%%%%%%%%%%%%%%%%%%%%%%%%%%%%%%
% %%%%%%%%%%%%%%%%%%%%%%%%%%%%%%%%%%%%%%%%%%%%%%%%%

\section{Examples}
\label{sec:examples}

This section contains examples of the theory we developed in the earlier sections.
We work under the following setup.

\begin{setup}
\label{setup:A4-examples}
Let $k$ be an algebraically closed field. 
The following examples come from the bounded derived category of the path algebra $kQ$, where $Q$ is the linearly-oriented Dynkin quiver $1\to2\to3\to4$. 
We denote by $\CC$ the bounded derived category $D^{b}(\lmod{kQ})$ 
of the category of finite-dimensional left $kQ$-modules with shift functor $[1]$. 
We will refer to the following portion of the Auslander-Reiten quiver of $\CC$ in these examples.
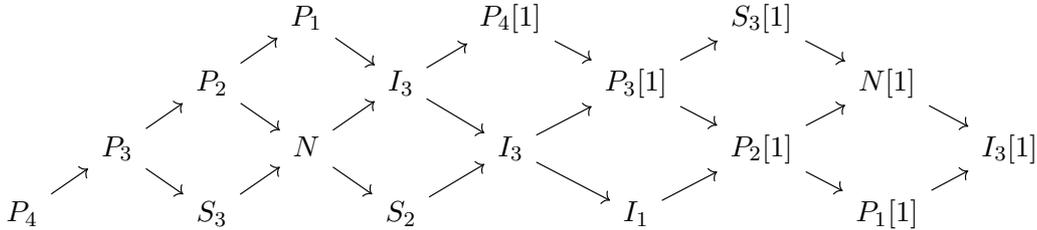
\begin{figure}[H]
\centering
\[
\begin{tikzcd}[column sep=0.5cm, row sep=0.2cm, scale cd=1]
&&& P_{1} 
    \arrow{dr}
&&P_{4}[1] 
    \arrow{dr}
&&S_{3}[1]
    \arrow{dr}
&&
\\
&& P_{2}
    \arrow{ur}
    \arrow{dr}
&& I_{2}
    \arrow{dr} 
    \arrow{ur}
&& P_{3}[1] 
    \arrow{ur} 
    \arrow{dr}
&& N[1]
    \arrow{dr}
&
\\
& P_{3} 
    \arrow{ur} 
    \arrow{dr}
&& N 
    \arrow{ur}
    \arrow{dr} 
&& I_{3} 
    \arrow{dr}
    \arrow{ur}
&& P_{2}[1]
    \arrow{dr}
    \arrow{ur}
&&I_{2}[1]  
\\
P_{4} 
    \arrow{ur} 
&& S_{3} 
    \arrow{ur}
&& S_{2}
    \arrow{ur} 
&& I_{1} 
    \arrow{ur}
&& P_{1}[1]
    \arrow{ur}
&
\end{tikzcd}
\]
\caption{Segment of the Auslander-Reiten quiver of $\CC$.}
\label{fig:AR-quiver-segment}
\end{figure}

Lastly, we let $\A$ denote some full subcategory of $\CC$ that is extension-closed, closed under direct summands, and contains the indecomposable objects depicted in Figure~\ref{fig:AR-quiver-segment}. 
It follows that $\A$ admits an extriangulated structure, say $(\BE,\fs)$, which is inherited from the triangulated structure of $\CC$. 
Further, $(\A,\BE,\fs)$ is an artin $k$-linear extriangulated category. 
Unless  $\A=\CC$, then $\A$ will not be closed under the shift functor $[1]$ of $\CC$ and so will not be a triangulated category. 
\end{setup}

In the examples in this section, when we say `$(\restr{\BE}{\F(\Phi)},\restr{\fs}{\F(\Phi)})$-subobject' we mean $(\restr{\BE}{\F(\Phi)},\restr{\fs}{\F(\Phi)})$-subobject up to isomorphism. 
The examples we give are not exact categories and are thus not covered under the theory developed in \cite{MS} or that in \cite{Sa}.

\begin{example}
\label{example:strong-proj-ESS}
Let 
$\Phi = \{ \Phi_{i} \}_{i=1}^{3}$ 
where 
$\Phi_{1} \deff S_{2}$, 
$\Phi_{2} \deff P_{3}$ and  
$\Phi_{3} \deff S_{3}[1]$. 
It follows that 
$\F(\Phi) = \add\set{\Phi_{1},\Phi_{2},\Phi_{3},P_{2}}$. 
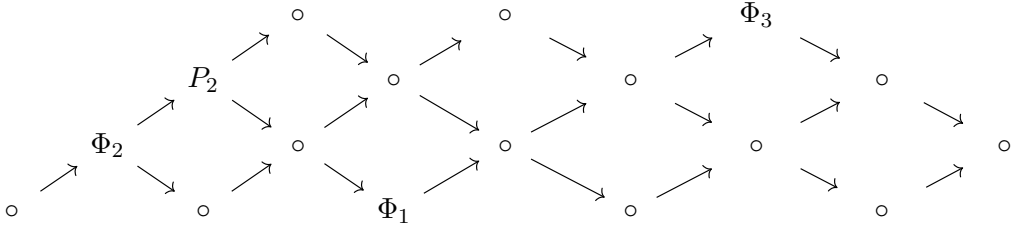
\[
\begin{tikzcd}[column sep=0.5cm, row sep=0.2cm, scale cd=1]
&&& \overlap{\circ}{\phantom{P_{1}}}
    \arrow{dr}
&&  \overlap{\circ}{\phantom{P_{4}[1]}}
    \arrow{dr}
&& \overlap{\Phi_{3}}{\phantom{S_{3}[1]}}
    \arrow{dr}
&&
\\
&& \overlap{P_{2}}{\phantom{P_{2}}}
    \arrow{ur}
    \arrow{dr}
&& \overlap{\circ}{\phantom{I_{3}}}
    \arrow{dr} 
    \arrow{ur}
&& \overlap{\circ}{\phantom{P_{3}[1]}}
    \arrow{ur} 
    \arrow{dr}
&&  \overlap{\circ}{\phantom{N[1]}}
    \arrow{dr}
&
\\
& \overlap{\Phi_{2}}{\phantom{P_{3}}}
    \arrow{ur} 
    \arrow{dr}
&& \overlap{\circ}{\phantom{N}}
    \arrow{ur}
    \arrow{dr} 
&& \overlap{\circ}{\phantom{I_{3}}}
    \arrow{dr}
    \arrow{ur}
&& \overlap{\circ}{\phantom{P_{2}[1]}}
    \arrow{dr}
    \arrow{ur}
&& \overlap{\circ}{\phantom{I_{3}[1]}}
\\
\overlap{\circ}{\phantom{P_{4}}}
    \arrow{ur} 
&& \overlap{\circ}{\phantom{S_{3}}}
    \arrow{ur}
&& \Phi_{1}
    \arrow{ur} 
&& \overlap{\circ}{\phantom{I_{1}}}
    \arrow{ur}
&& \overlap{\circ}{\phantom{P_{1}[1]}}
    \arrow{ur}
&
\end{tikzcd}
\]  
Set $Q = \{ Q_{i} \}_{i=1}^{3}$ 
where $Q_{1} \deff P_{2}, Q_{2} \deff \Phi_{2}, Q_{3} \deff \Phi_{3}$. 
Notice that the only indecomposable extriangle in $(\F(\Phi), \restr{\BE}{\F(\Phi)},\restr{\fs}{\F(\Phi)})$ is 
$\begin{tikzcd}[column sep=0.5cm,cramped]
P_{3} \arrow{r}{k_{1}}& P_{2} \arrow{r}{q_{1}}& S_{2} \arrow[dashed]{r}{\eta_{1}}&{},
\end{tikzcd}$
where $q_{1}$ is a non-zero morphism. 
Let $\eta_{i} = {}_{0}0_{\Phi_{i}}$ so that $q_{i} = \id{\Phi_{i}}$ for $i=2,3$.
Using this and Figure~\ref{fig:AR-quiver-segment}, one easily checks that $(\Phi, Q)$ is a minimal projective $\BE$-stratifying system.

We claim that $(\F(\Phi), \restr{\BE}{\F(\Phi)},\restr{\fs}{\F(\Phi)})$ is a Jordan-H{\"{o}}lder extriangulated category that is not an exact category. 
It is straightforward to see that 
$\A(Q_{i}, -)$ is left exact on $(\F(\Phi), \restr{\BE}{\F(\Phi)},\restr{\fs}{\F(\Phi)})$ for all $1\leq i \leq 3$ by checking left exactness on the (up to isomorphism) single indecomposable extriangle.
Therefore, by \cref{cor:JHP-for-F-Theta}, the category $(\F(\Phi), \restr{\BE}{\F(\Phi)},\restr{\fs}{\F(\Phi)})$ has the Jordan-H{\"{o}}lder property with $(\restr{\BE}{\F(\Phi)},\restr{\fs}{\F(\Phi)})$-simple objects $\Phi_{1}, \Phi_{2}, \Phi_{3}$. To see that $(\F(\Phi), \restr{\BE}{\F(\Phi)},\restr{\fs}{\F(\Phi)})$ is not exact, observe that the deflation $q_{1}\colon P_{2} \to S_{2}$ is not an epimorphism. Indeed, let $y\colon S_{2} \to S_{3}[1]$ be a non-zero map, then $yq_{1}=0$.
\end{example}

\begin{remark}
\label{rem:counter-to-EnomotoSaito-Q6-3}
In a preliminary version of \cite{EnomotoSaito}, it was asked if a skeletally small extriangulated category is exact if and only if it has a reduced Grothendieck monoid.\footnote{See Question 6.3 in \url{https://arxiv.org/abs/2208.02928v1}.} 
The category $(\F(\Phi), \restr{\BE}{\F(\Phi)},\restr{\fs}{\F(\Phi)})$ in \cref{example:strong-proj-ESS} above provides a negative answer. 
Indeed, $0$ is $(\restr{\BE}{\F(\Phi)},\restr{\fs}{\F(\Phi)})$-simple-like, so 
the Grothendieck monoid 
$\M(\F(\Phi), \restr{\BE}{\F(\Phi)},\restr{\fs}{\F(\Phi)})$ is reduced by 
\cref{prop:consequences-of-zero-simple-like}\ref{item:M-reduced}, but $(\F(\Phi), \restr{\BE}{\F(\Phi)},\restr{\fs}{\F(\Phi)})$ is not an exact category as it has a deflation that is not an epimorphism. 
The authors thank Haruhisa Enomoto and Shunya Saito for bringing this to their attention.
\end{remark}

\begin{example} \label{example:zero projective}
Let 
$\Phi = \{ \Phi_{i} \}_{i=1}^{3}$ 
where 
$\Phi_{1} \deff P_{2}[1]$, 
$\Phi_{2} \deff S_{2}$ and  
$\Phi_{3} \deff P_{3}$. 
It follows that 
$\F(\Phi) = \add\set{\Phi_{1},\Phi_{2},\Phi_{3}, P_{2}, P_{3}[1]}$. 
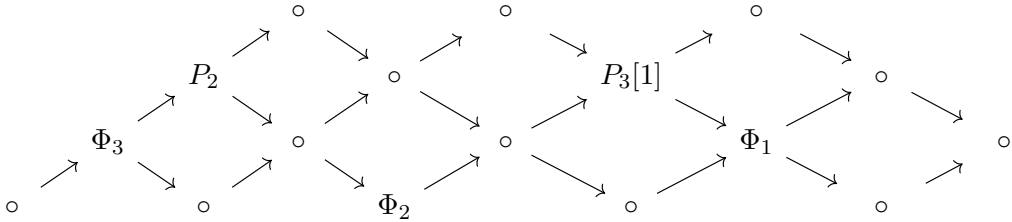
\[
\begin{tikzcd}[column sep=0.5cm, row sep=0.2cm, scale cd=1]
&&& \overlap{\circ}{\phantom{P_{1}}}
    \arrow{dr}
&&  
    \overlap{\circ}{\phantom{P_{4}[1]}}
    \arrow{dr}
&& \overlap{\circ}{\phantom{S_{3}[1]}} 
    \arrow{dr}
&&
\\
&& P_{2}
    \arrow{ur}
    \arrow{dr}
&& \overlap{\circ}{\phantom{I_{3}}}
    %\circ %I_{3}
    \arrow{dr} 
    \arrow{ur}
&& P_{3}[1]
    % \circ %P_{3}[1] 
    \arrow{ur} 
    \arrow{dr}
&& \overlap{\circ}{\phantom{\Phi_{1}}}
    \arrow{dr}
&
\\
& \Phi_{3} 
    \arrow{ur} 
    \arrow{dr}
&& \overlap{\circ}{\phantom{N}}
    % \circ %N 
    \arrow{ur}
    \arrow{dr} 
&& \overlap{\circ}{\phantom{I_{3}}}
    % \circ %I_{3} 
    \arrow{dr}
    \arrow{ur}
&& \Phi_{1}
    % \circ %P_{2}[1]
    \arrow{dr}
    \arrow{ur}
&& \overlap{\circ}{\phantom{I_{3}[1]}}
    % \circ %I_{3}[1]  
\\
\overlap{\circ}{\phantom{P_{4}}}
% \circ %P_{4} 
    \arrow{ur} 
&& \overlap{\circ}{\phantom{S_{3}}}
    % \circ %S_{3} 
    \arrow{ur}
&& \Phi_{2}
    \arrow{ur} 
&& \overlap{\circ}{\phantom{I_{1}}}
% \circ %I_{1} 
    \arrow{ur}
&& \overlap{\circ}{\phantom{P_{1}[1]}}
% \circ %P_{1}[1]
    \arrow{ur}
&
\end{tikzcd} 
\] 
Set $Q = \{ Q_{i} \}_{i=1}^{3}$ 
where 
$Q_{1} \deff 0$, $Q_{2}\deff P_{2}$, $Q_{3}\deff \Phi_{3}$. 
We have the extriangles 
$\begin{tikzcd}[column sep=0.5cm,cramped]
P_{2} \arrow{r}& 0 \arrow{r}{q_{1}}& \Phi_{1} \arrow[dashed]{r}{\eta_{1}}& {}
\end{tikzcd}$
and 
$\begin{tikzcd}[column sep=0.5cm,cramped]
\Phi_{3} \arrow{r}& P_{2} \arrow{r}{q_{2}}& \Phi_{2} \arrow[dashed]{r}{\eta_{2}}& {,}
\end{tikzcd}$
and we put $\eta_{3} = {}_{0}0_{\Phi_{3}}$.
It is straightforward to verify that $(\Phi, Q)$ is a minimal projective $\BE$-stratifying system.

We claim that 
$(\F(\Phi), \restr{\BE}{\F(\Phi)},\restr{\fs}{\F(\Phi)})$ is not Jordan-H{\"{o}}lder and also not an exact category. 
Both claims follow from the 
existence of the extriangle 
$\begin{tikzcd}[column sep=0.5cm,cramped]
P_{3} \arrow{r}& 0\arrow{r} & P_{3}[1] \arrow[dashed]{r}&{}
\end{tikzcd}$
in $(\F(\Phi), \restr{\BE}{\F(\Phi)},\restr{\fs}{\F(\Phi)})$. 
This extriangle is clearly not a short exact sequence and thus $(\F(\Phi), \restr{\BE}{\F(\Phi)},\restr{\fs}{\F(\Phi)})$ is not an exact category. 
Furthermore, as we discussed in \cref{sec:Grothendieck-monoid-and-JH-property}, this extriangle allows for infinitely many filtrations of $P_{3}$ of different lengths---an obstruction to the Jordan-H{\"{o}}lder property. Thus, by \cref{cor:JHP-for-F-Theta}, as 
$(\Phi, Q)$ is a minimal projective $\BE$-stratifying system, we must have that $\A(Q_{i},-)$ is not left exact on $(\F(\Phi), \restr{\BE}{\F(\Phi)},\restr{\fs}{\F(\Phi)})$ for some $i$. Indeed, $\A(Q_{3},-)$ is not left exact as seen using the extriangle $\begin{tikzcd}[column sep=0.5cm,cramped]
P_{3} \arrow{r}& 0\arrow{r} & P_{3}[1] \arrow[dashed]{r}&{.}
\end{tikzcd}$
\end{example}

In our final example, we show that the left exactness condition of \cref{cor:JHP-for-F-Theta} is not necessary for $(\F(\Phi), \restr{\BE}{\F(\Phi)},\restr{\fs}{\F(\Phi)})$ to be Jordan-H{\"{o}}lder.

\begin{example}
Let 
$\Phi = \{ \Phi_{i} \}_{i=1}^{3}$ 
where 
$\Phi_{1} \deff N[1]$, 
$\Phi_{2} \deff S_{2}$ and  
$\Phi_{3} \deff P_{3}$. 
It follows that 
$\F(\Phi) = \add\set{\Phi_{1},\Phi_{2},\Phi_{3}, P_{2}, S_{3}[1]}$. 
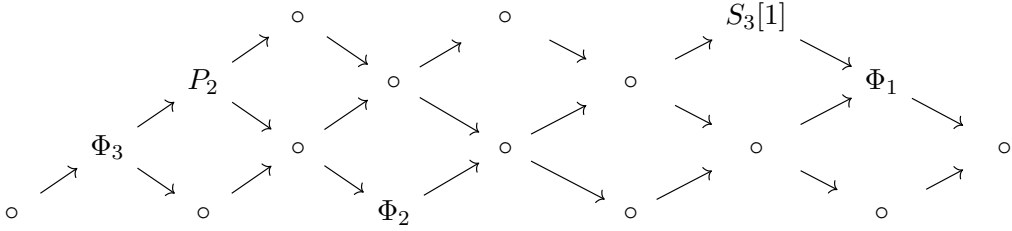
\[
\begin{tikzcd}[column sep=0.5cm, row sep=0.2cm, scale cd=1]
&&& \overlap{\circ}{\phantom{P_{1}}}
    %\circ %P_{1} 
    \arrow{dr}
&&  %\phantom{P_{4}[1]}
    % P_{4}[1]
    % {\makebox[\widthof{$P_{4}[1]$}/2-\widthof{$\circ$}/2][l]{$P_{4}[1]$}\circ}
    \overlap{\circ}{\phantom{P_{4}[1]}}
    \arrow{dr}
&& S_{3}[1]
    \arrow{dr}
&&
\\
&& P_{2}
    \arrow{ur}
    \arrow{dr}
&& \overlap{\circ}{\phantom{I_{3}}}
    %\circ %I_{3}
    \arrow{dr} 
    \arrow{ur}
&& \overlap{\circ}{\phantom{P_{3}[1]}}
    % \circ %P_{3}[1] 
    \arrow{ur} 
    \arrow{dr}
&& \Phi_{1}
    \arrow{dr}
&
\\
& \Phi_{3} 
    \arrow{ur} 
    \arrow{dr}
&& \overlap{\circ}{\phantom{N}}
    % \circ %N 
    \arrow{ur}
    \arrow{dr} 
&& \overlap{\circ}{\phantom{I_{3}}}
    % \circ %I_{3} 
    \arrow{dr}
    \arrow{ur}
&& \overlap{\circ}{\phantom{P_{2}[1]}}
    % \circ %P_{2}[1]
    \arrow{dr}
    \arrow{ur}
&& \overlap{\circ}{\phantom{I_{3}[1]}}
    % \circ %I_{3}[1]  
\\
\overlap{\circ}{\phantom{P_{4}}}
% \circ %P_{4} 
    \arrow{ur} 
&& \overlap{\circ}{\phantom{S_{3}}}
    % \circ %S_{3} 
    \arrow{ur}
&& \Phi_{2}
    \arrow{ur} 
&& \overlap{\circ}{\phantom{I_{1}}}
% \circ %I_{1} 
    \arrow{ur}
&& \overlap{\circ}{\phantom{P_{1}[1]}}
% \circ %P_{1}[1]
    \arrow{ur}
&
\end{tikzcd}
\] 
Set $Q = \{ Q_{i} \}_{i=1}^{3}$ 
where 
$
Q_{1} \deff S_{3}[1],
Q_{2} \deff P_{2},
Q_{3} \deff \Phi_{3}
$. 
We have the extriangles
\begin{align}
&\begin{tikzcd}[ampersand replacement=\&]
\Phi_{2} \arrow[r, "{k_{1}}"] 
	\& {S_{3}[1]} \arrow[r, "{q_{1}}"] 
	\& \Phi_{1} \arrow[r, "{\eta_{1}}", dashed] 
	\& {,} 
\end{tikzcd}\label{eqn:example-5.5-eta1}
\\
&\begin{tikzcd}[ampersand replacement=\&]
\Phi_{3} \arrow[r, "{k_{2}}"]  \& P_{2} \arrow[r, "{q_{2}}"]       \& \Phi_{2} \arrow[r, "{\eta_{2}}", dashed]     \& {,}
\end{tikzcd}\label{eqn:example-5.5-eta2}
\\
&\begin{tikzcd}[ampersand replacement=\&]
0 \arrow{r}\& \Phi_{3} \arrow{r}{\id{\Phi_{3}}}\& \Phi_{3} \arrow[dashed]{r}{\eta_{3}}\& {,}
\end{tikzcd}\label{eqn:example-5.5-eta3}
\end{align}
where $k_{1},k_{2},q_{1},q_{2},q_{3}$ are all right minimal.

We claim that 
$(\Phi, Q)$ is a minimal projective $\BE$-stratifying system, 
that $(\F(\Phi), \restr{\BE}{\F(\Phi)},\restr{\fs}{\F(\Phi)})$ is Jordan-H{\"{o}}lder but not exact, and 
that $\A(Q_{2}, -)$ is not left exact on $(\F(\Phi), \restr{\BE}{\F(\Phi)},\restr{\fs}{\F(\Phi)})$. 
As in previous examples, one may verify \ref{S1}, \ref{S2}, \ref{PS1} hold using Figure~\ref{fig:AR-quiver-segment}. 
The only indecomposable extriangles (up to isomorphism) in $(\F(\Phi), \restr{\BE}{\F(\Phi)},\restr{\fs}{\F(\Phi)})$ are the extriangles \eqref{eqn:example-5.5-eta1} and \eqref{eqn:example-5.5-eta2}.
Although each $q_{i}$ is non-zero for $i=1,2,3$, 
the hypotheses of \cref{cor:JHP-for-F-Theta} are not fully satisfied
since the functor
$
\A(Q_{2},-) = \A(P_{2},-)
$ 
is not left exact on $(\F(\Phi), \restr{\BE}{\F(\Phi)},\restr{\fs}{\F(\Phi)})$. 
Nevertheless, $(\F(\Phi), \restr{\BE}{\F(\Phi)},\restr{\fs}{\F(\Phi)})$ enjoys the Jordan-H{\"{o}}lder property, which we can verify. Since the only extriangles in $(\F(\Phi), \restr{\BE}{\F(\Phi)},\restr{\fs}{\F(\Phi)})$ with decomposable middle terms are split extriangles, we only have to concern ourselves with indecomposable objects. 
Of these, only $P_{2}$ and $S_{3}[1]$ are not $(\restr{\BE}{\F(\Phi)},\restr{\fs}{\F(\Phi)})$-simple, but they each have a unique $(\restr{\BE}{\F(\Phi)},\restr{\fs}{\F(\Phi)})$-subobject. In both cases, the $(\restr{\BE}{\F(\Phi)},\restr{\fs}{\F(\Phi)})$-subobject is $(\restr{\BE}{\F(\Phi)},\restr{\fs}{\F(\Phi)})$-simple and the Jordan-H{\"{o}}lder property follows.

Lastly, we note that the deflation $q_{2}$ is not an epimorphism, because $k_{1}q_{2} = 0$ but $k_{1}\neq 0$. Hence, $(\F(\Phi), \restr{\BE}{\F(\Phi)},\restr{\fs}{\F(\Phi)})$ is not an exact category. 
\end{example}

%%%%%%%%%%%%%%%%%%%%%%%%%%%%%%%%%%%%%%
%%%%%%%%%%%%%%%%%%%%%%%%%%%%%%%%%%%%%%
% Acknowledgements
{\setstretch{1}
\begin{acknowledgements}
We thank Mikhail Gorsky for interesting discussions. 
We are very grateful to Carlo Klapproth for bringing \cref{prop:A-has-simple-then-zero-simple} to their attention, which simplified several results including \cref{thm:JH-iff-free-monoid-iff-free-group}.
We are also grateful to Haruhisa Enomoto and Shunya Saito for email communications, which led to \cref{rem:counter-to-EnomotoSaito-Q6-3}.
We thank an anonymous referee whose comments led to several improvements to the paper, especially to a stronger \cref{cor:JHP-for-F-Theta}. We also thank a second referee for a careful reading of the paper and spotting several typos.

The first author was supported by NSERC of Canada.
The third author gratefully acknowledges support from: the Danish National Research Foundation (grant DNRF156); the Independent Research Fund Denmark (grant 1026-00050B); the Aarhus University Research Foundation (grant AUFF-F-2020-7-16); the Engineering and Physical Sciences Research Council (grant EP/P016014/1); and the London Mathematical Society with support from Heilbronn Institute for Mathematical Research (grant ECF-1920-57).
\end{acknowledgements}}

%%%%%%%%%%%%%%%%%%%%%%%%%%%%%%%%%%%%%%
%%%%%%%%%%%%%%%%%%%%%%%%%%%%%%%%%%%%%%
% References
{\setstretch{1}
%\bibliography{references}
%\bibliographystyle{mybstwithlabels}

}
\end{document}